\documentclass[a4paper,12pt]{article}
\pdfoutput=1
   \usepackage[top=2.5cm,bottom=2.5cm,left=2.5cm,right=2.5cm]{geometry}
\usepackage{cite,amsmath,amssymb}
\usepackage[margin=1cm,%
            font=small,%
            format=hang,%
            labelsep=period,%
             labelfont=bf]{caption}
\usepackage{graphicx}
\usepackage{caption}
\usepackage{float}
\usepackage{amsmath}
\usepackage{cases}
\usepackage{mathrsfs}
\usepackage{amsfonts}
\usepackage{color}
\usepackage{setspace}
\usepackage{cite}
\usepackage{amsmath}
\usepackage{array}
\usepackage{booktabs}

\usepackage{graphicx}
\usepackage{subfigure}
\usepackage{CJK}
\usepackage{caption}
\usepackage{float}
\usepackage{amsmath}
\usepackage{cases}
\usepackage{mathrsfs}
\usepackage{amsfonts}
\usepackage{times}
\usepackage{mathptmx}
\usepackage{cite}
\usepackage{setspace}
\usepackage{times}
\usepackage{mathptmx}
\usepackage{geometry}

\renewcommand{\figurename}

\newtheorem{theorem}{Theorem}[section]
\newtheorem{lemma}{Lemma}[section]

\renewcommand{\figurename}

\makeatletter
\def\thanks#1{\protected@xdef\@thanks{\@thanks
        \protect\footnotetext{#1}}}
\makeatother

\begin{document}

\title{\textbf{Arithmetic-Geometric spectral radii of Unicyclic graphs}
\vspace*{0.3cm}}
\author{ \bf  Ruiling Zheng, \, Xian$^{^{\textbf{,}}}$an Jin\\
\small\em School of Mathematical Sciences, Xiamen University, Xiamen  361005, P. R. China\\ \small\em E-mail:
rlzheng2017@163.com, xajin@xmu.edu.cn.\vspace*{0.3cm}}
 \date{}
\maketitle\thispagestyle{empty}
\vspace*{-1cm}

\begin{abstract}
Let $d_{v_{i}}$ be the degree of the vertex $v_{i}$ of $G$. The arithmetic-geometric matrix $A_{ag}(G)$ of a graph $G$ is a square matrix, where the $(i,j)$-entry is equal
to $
\displaystyle \frac{d_{v_{i}}+d_{v_{j}}}{2\sqrt{d_{v_{i}}d_{v_{j}}}}$ if the vertices $v_{i}$ and $v_{j}$ are adjacent, and 0
otherwise. The
arithmetic-geometric spectral radius of $G$, denoted by $\rho_{ag}(G)$, is the largest eigenvalue of the
arithmetic-geometric matrix $A_{ag}(G)$. In this paper, the unicyclic graphs of order $n\geq5$ with the smallest and first four largest arithmetic-geometric spectral radii are determined.

\medskip

{\bf Keywords:} \  Arithmetic-geometric matrix; Arithmetic-geometric spectral radius; Unicyclic graphs; Smallest and first four largest
\end{abstract}

\vspace*{0.35cm}
\baselineskip=0.30in

\section{Introduction}
\label{intro}

Throughout this paper we consider finite, undirected, simple and connected graphs. Let $G$ be a graph with vertex set
$V(G)=\{v_{1},v_{2},\ldots,v_{n}\}$ and edge set $E(G)$. An edge $e\in E(G)$ with end vertices $v_{i}$ and $v_{j}$
is denoted by
$v_{i}v_{j}$. For $i=1,2,\ldots,n$, we denote by $d_{v_{i}}$ the degree of the vertex $v_{i}$ in $G$. For a vertex $v$ of a graph $G$, let $N_{G}(v)=\{u\in V(G):uv\in E(G)\}$. As usual, let
$S_{n}$ and
$C_{n}$ be the star and the cycle of order $n\geq3$. And $S_{n}+e$ with $n\geq3$ denote the unicyclic graph obtained from $S_{n}$ by adding an edge, see Fig. \ref{fig1}. The spectral radius of a complex matrix $M$ is the largest value among the modules
of all eigenvalues of $M$. If $M\geq0$ and it is irreducible and real symmetric, then its spectral radius is
exactly the largest eigenvalue and we denote it by $\displaystyle \rho(M)$. The spectral radius of the adjacency
matrix
$A(G)=(a_{i,j})$ of $G$ is referred as the spectral radius of $G$.

In 1975, a so-called Randi\'{c} index was proposed \cite{1}, it is
defined as

\begin{center}
$ \displaystyle R(G)=\sum\limits_{v_{i}v_{j}\in E(G)}\frac{1}{\sqrt{d_{v_{i}}d_{v_{j}}}}.$
\end{center}

\noindent The Randi\'{c} matrix $R(G)=(r_{i,j})$ of $G$, where $\displaystyle
r_{i,j}=\frac{1}{\sqrt{d_{v_{i}}d_{v_{j}}}}$ if $v_{i}v_{j}\in E(G)$ and 0 otherwise, seems to be
first time used in 2005 by Rodr\'{\i}guez, who referred to it as the ``weighted adjacency matrix" \cite{555} and
the ``degree adjacency matrix" \cite{666}. The Randi\'{c} spectral radius was studied in \cite{77,999,2,888} and
the references cited therein.

In 1998, Estrada, Torres, Rodr\'{\i}guez and Gutman \cite{8} introduced the atom-bond connectivity (ABC) index,
defined as

\begin{center}
$\displaystyle ABC(G)=\sum\limits_{v_{i}v_{j}\in E(G)}\sqrt{\frac{d_{v_{i}}+d_{v_{j}}-2}{d_{v_{i}}d_{v_{j}}}}.$
\end{center}

\noindent The $ABC$ matrix of a graph $G$, put forward by Estrada \cite{10}, is defined to be $\Omega(G)=(\omega_{i,j})$,
where
$\displaystyle \omega_{i,j}=\sqrt{\frac{d_{v_{i}}+d_{v_{j}}-2}{d_{v_{i}}d_{v_{j}}}}$ if $v_{i}v_{j}\in E(G)$ and $0$ otherwise. Li
and Wang \cite{13} proved that for an unicyclic graph $G$ of order $n\geq4$,

\begin{center}
$\displaystyle \sqrt{2}=\rho(\Omega(C_{n}))\leq\rho(\Omega(G))\leq\rho(\Omega(S_{n}+e)),$
\end{center}

\noindent with equality if and only if $T\cong C_{n}$ for the lower bound, and if and only if $T\cong S_{n}+e$ for
the upper bound,
which was conjectured in \cite{14}. Recently, Yuan et al. \cite{462} considered the unicyclic graphs with the first
four largest $ABC$ spectral radii.

In 2015, Shegehall and Kanabur \cite{15} proposed the arithmetic-geometric index of $G$, which is defined as

\begin{center}
$ \displaystyle AG(G)=\sum\limits_{v_{i}v_{j}\in E(G)}\frac{d_{v_{i}}+d_{v_{j}}}{2\sqrt{d_{v_{i}}d_{v_{j}}}}.  $
\end{center}

\noindent   Zheng et al. \cite{16} considered the arithmetic-geometric matrix ($AG$ matrix) of $G$, denoted by
$A_{ag}(G)=(h_{i,j})$,
where $\displaystyle h_{i,j}=\frac{d_{v_{i}}+d_{v_{j}}}{2\sqrt{d_{v_{i}}d_{v_{j}}}}$ if $v_{i}v_{j}\in E(G)$ and 0 otherwise. We
denote the $\rho(A_{ag}(G))$  by $\rho_{ag}(G)$ and call it
the arithmetic-geometric spectral radius ($AG$ spectral radius) of $G$. And for a graph $G$, let $\Phi(G,\rho)$ be the characteristic polynomial of the $AG$ matrix of $G$. Zheng et al. \cite{16} obtained two lower
bounds for it in terms of the order, maximum degree and first Zagreb index $M_{1}(G)=\sum\limits_{v_{i}v_{j}\in
E}(d_{v_{i}}+d_{v_{j}})$ \cite{M1} of $G$. And they also proved that for a graph $G$ of order $n\geq2$ and size $m$,

\begin{center}
$\displaystyle\rho_{ag}(G)\leq\frac{1}{2}\Big(\sqrt{n-1}+\frac{1}{\sqrt{n-1}}\Big)\sqrt{2m-n+1},$
\end{center}
 where the equality holds if and only if $G\cong S_{n}$. Guo and Gao \cite{17} offered several other bounds for the $AG$ spectral radius in terms of the order and size, the maximum, minimum degree and the spectral radius
 $\displaystyle \rho(G)$ of $G$ as well.

In this paper, we determine the unicyclic graphs of order $n\geq5$ with the smallest and first four largest
arithmetic-geometric spectral radii.

\section{Lemmas}
\label{lem}

We need the following lemmas.

\begin{lemma}\label{lemma 1}\cite{19}
Let $M$ ba a nonnegative matrix of order $n$. Let $\textbf{x}$ be a positive column vector of dimension
$n$, i.e.
every entry of $\textbf{x}$ is positive. If $k>0$ such that $M\textbf{x}\leq k\textbf{x}$, then $\rho(M)\leq k$.
\end{lemma}

\begin{lemma}\label{lemma 2}\cite{19}
Let $B=(b_{ij})$ and $D=(d_{ij})$ be two nonnegative matrices of order $n$. If $B\geq D$, i.e., $b_{ij} \geq
d_{ij}\geq0$
for all $i,j$, then $\rho(B) \geq\rho(D)$.
\end{lemma}

\begin{lemma}\label{lemma 3}\cite{22}
If $G$ is an unicyclic graph of order $n\geq3$, then
\begin{center}
$\displaystyle 2=\rho(A(C_{n}))\leq\rho(A(G))\leq\rho(A(S_{n}+e)),$
\end{center}

\noindent the lower bound is attained if and only if $G\cong C_{n}$, the upper bound is obtained if and only if
$G\cong S_{n}+e$.
\end{lemma}

\section{Main Result}
\label{sec:2}

In this section, we study the $AG$ spectral radius of unicyclic graphs.

\begin{lemma}\label{lemma 4} Let $G$ be a graph with $m$ edges without isolated vertices. Let $\displaystyle V(G)=\{v_{1},...,v_{n}\}$ and $\displaystyle X=(\sqrt{d_{v_{1}}},...,\sqrt{d_{v_{n}}})^{T}$. Then for any vertex $\displaystyle v_{i}\in V(G)$, we have
\begin{center}
$\displaystyle (A_{ag}(G)X)_{v_{i}}\leq \frac{d_{v_{i}}^{2}+2m-n+1}{2\sqrt{d_{v_{i}}}}$.
\end{center}
\end{lemma}

\noindent \textbf{Proof}  Note that for any vertex $\displaystyle  v_{i}\in V(G)$, we have

\begin{eqnarray*}
(A_{ag}(G)X)_{v_{i}}&=&\sum_{v_{j}\in N_{G}(v_{i})}\frac{d_{v_{i}}+d_{v_{j}}}{2\sqrt{d_{v_{i}}d_{v_{j}}}}\cdot
 \sqrt{d_{v_{j}}}=\frac{1}{2}\sum_{v_{j}\in N_{G}(v_{i})}(\sqrt{d_{v_{i}}}+\frac{d_{v_{j}}}{\sqrt{d_{v_{i}}}})\\
 &=&\frac{d_{v_{i}}\sqrt{d_{v_{i}}}}{2}+\frac{\sum_{v_{j}\in N_{G}(v_{i})}d_{v_{j}}}{2\sqrt{d_{v_{i}}}}.
\end{eqnarray*}

\noindent Since

\begin{eqnarray*}
\sum_{v_{j}\in N_{G}(v_{i})}d_{v_{j}}&=&2m-d_{v_{i}}-\sum_{v_{k}\notin (N_{G}(v_{i})\cup\{v_{i}\})}d_{v_{k}}\\
&\leq&2m-d_{v_{i}}-\sum_{v_{k}\notin N_{G}(v_{i})\cup\{v_{i}\}}1\\
&=& 2m-d_{v_{i}}-(n-1-d_{v_{i}})\\
&=& 2m-n+1.
\end{eqnarray*}

\noindent Thus we obtain

\begin{center}
$\displaystyle (A_{ag}(G)X)_{v_{i}}\leq \frac{d_{v_{i}}^{2}+2m-n+1}{2\sqrt{d_{v_{i}}}}$.
\end{center}\hfill$\Box$

\begin{lemma}\label{lemma 7} Let $G$ be an unicyclic graph of order $n\geq8$ with maximum degree $\Delta$. Then

$$ \rho_{ag}(G)<\left\{
\begin{aligned}
&\frac{n-1.45}{2}, \ \ if\  8\leq n\leq 15\ and \ \Delta\leq n-4;\\
&\frac{n-1.65}{2}, \ \ if \  16\leq n\leq 21\  and \ \Delta\leq n-3;\\
&\frac{n-1.75}{2}, \ \ if \  22\leq n \ and \ \Delta\leq n-3.
\end{aligned}
\right.
$$
\end{lemma}

\noindent \textbf{Proof} We only prove the case of $n\leq 22$ and $\Delta\leq n-3$, and the other two cases can be proved similarly.

Let $\displaystyle X=(\sqrt{d_{v_{1}}},...,\sqrt{d_{v_{n}}})^{T}$. If $d_{v_{i}}=1$, let $v_{j}$ is the only neighbor of $v_{i}$, then

\begin{center}
$\displaystyle (A_{ag}(G)X)_{v_{i}}=\frac{1+d_{v_{j}}}{2\sqrt{d_{v_{j}}}}\cdot\sqrt{d_{v_{j}}}=\frac{1+d_{v_{j}}}{2}\leq\frac{n-2}{2}<\frac{n-1.75}{2}$.
\end{center}

\noindent Assume that $2\leq d_{v_{i}}\leq n-3$. Then by Lemma \ref{lemma 4} with $m=n$, we have
\begin{center}
$\displaystyle (A_{ag}(G)X)_{v_{i}}\leq \frac{d_{v_{i}}^{2}+n+1}{2\sqrt{d_{v_{i}}}}$.
\end{center}

\noindent For the fixed $n$ and $\displaystyle 2\leq x\leq n-3$. Let

\begin{center}
$\displaystyle f(x)=x^{2}-(n-1.75)x+n+1$.
\end{center}

\noindent  Then we have
\begin{center}
$\displaystyle f(x)\leq max\{f(2), f(n-3)\}=max\{\frac{17-2n}{2}, \frac{19-n}{4}\}<0$
\end{center}

\noindent for $n\geq22$. That is

\begin{center}
$\displaystyle d_{v_{i}}^{2}-(n-1.75)d_{v_{i}}+n+1=f(d_{v_{i}})<0$
\end{center}

\noindent for $\displaystyle 2\leq d_{v_{i}}\leq n-3$ and $n\geq22$, which implies that

\begin{center}
$\displaystyle (A_{ag}(G)X)_{v_{i}}\leq \frac{d_{v_{i}}^{2}+n+1}{2\sqrt{d_{v_{i}}}}<\frac{(n-1.75)}{2}\sqrt{d_{v_{i}}}$.
\end{center}

\noindent Hence, by Lemma \ref{lemma 1}, we obtain $\displaystyle \rho_{ag}(G)<\frac{n-1.75}{2}$ for the unicyclic graph with $\displaystyle \Delta\leq n-3$ and $n\geq22$. \hfill$\Box$

\begin{figure}[H]
  \centering
    \includegraphics[width=10cm]{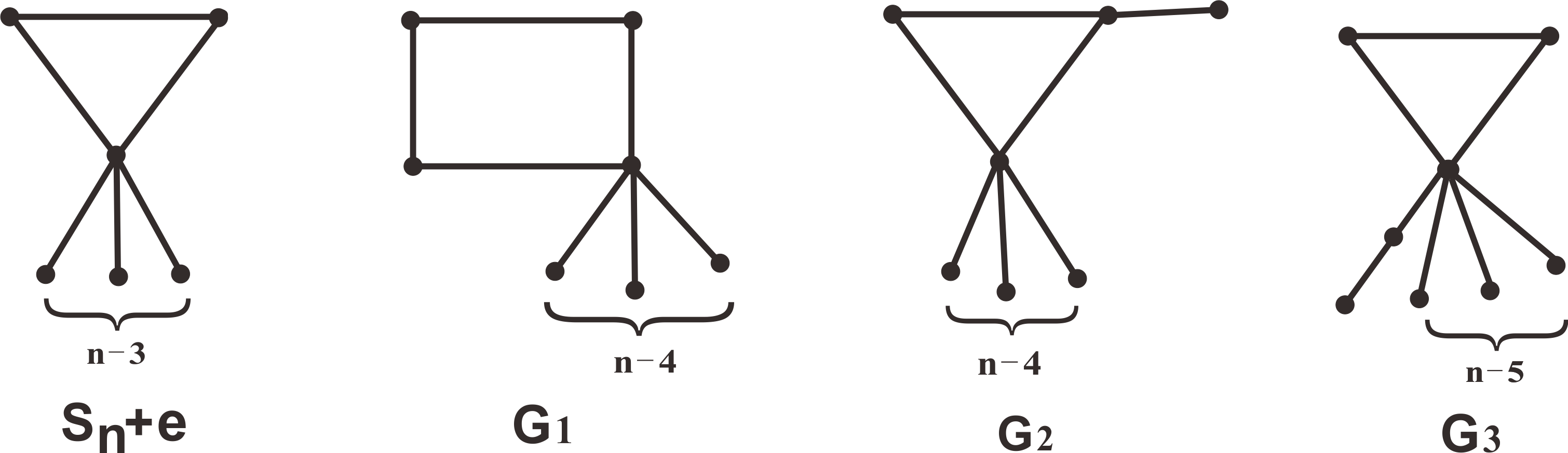}
  \caption{\small  The unicyclic graphs $S_{n}+e$ and $G_{i}$ with i={1,2,3}}
  \label{fig1}
\end{figure}

Let $G_{1}$ be the unicyclic graph of order $n\geq5$ obtained from a quadrangle by adding $n-4$ pendent edges incident to a common vertex, $G_{2}$ be the unicyclic graph of order $n\geq5$ obtained from $S_{n-1}+e$ by adding a pendent edge incident to a vertex of degree two, and $G_{3}$ be the unicyclic graph of order $n\geq5$ obtained from $S_{n-1}+e$ by adding a pendent edge incident to a vertex of degree one. See Fig. \ref{fig1}. We have the following lemma.

\begin{lemma}\label{lemma 6} For $n\geq8$, we have
\begin{itemize}
  \item[(i)] $\displaystyle \frac{n-1.45}{2} < \rho_{ag}(G_{3})<\rho_{ag}(G_{1})<\rho_{ag}(G_{2})< \frac{n-1}{2}$ if $8\leq n\leq15$;
  \item[(ii)] $\displaystyle \frac{n-1.65}{2} < \rho_{ag}(G_{3})<\rho_{ag}(G_{2})<\rho_{ag}(G_{1}) < \frac{n-1}{2}$ if $16 \leq n\leq 21$;
  \item[(iii)] $\displaystyle \frac{n-1.75}{2} < \rho_{ag}(G_{3})<\rho_{ag}(G_{2})<\rho_{ag}(G_{1}) < \frac{n-1}{2}$ if $n\geq 22$.
\end{itemize}
\end{lemma}
\noindent \textbf{Proof} For $8\leq n\leq 21$, the value of $\rho_{ag}(G_{1}),\rho_{ag}(G_{2}),\rho_{ag}(G_{3})$,  $\displaystyle \frac{n-1.65}{2}$, $\displaystyle \frac{n-1.45}{2}$ and $\displaystyle \frac{n-1}{2}$ are listed in Table \ref{table 1} of the Appendix 1. The results for (i) and (ii) can be obtained from Table \ref{table 1} directly.

Assume that $n\geq22$. Then (iii) follows from the following inequalities:
\begin{center}
$\displaystyle \frac{n-1.75}{2}<\rho_{ag}(G_{3})<\frac{n-1.6}{2}<\rho_{ag}(G_{2})<\frac{n-1.5}{2}<\rho_{ag}(G_{1})<\frac{n-1}{2}$.
\end{center}

By direct calculations, we obtain the arithmetic-geometric characteristic polynomial $\Phi_{i}(\rho,n)$ of $G_{i}$ $(i=1,2,3)$ as below:

\noindent $(1)$ $\displaystyle
\Phi_{1}(\rho,n)=\frac{\rho^{n-4}}{8(n-2)}[8(n-2)\rho^{4}-(2n^{3}-10n^{2}+34n-40)\rho^{2}+4(n-4)(n-1)^{2}]$;

\begin{spacing}{1.6}
\noindent $(2)$  $\displaystyle
\Phi_{2}(\rho,n)=\frac{\rho^{n-4}}{4(n-2)}\Big[4(n-2)\rho^{4}-\frac{6n^{3}-31n^{2}+115n-136}{6}\rho^{2}-\frac{5n^{2}+5n}{6}\rho+\frac{2n^{2}}{3}
+\frac{19(n-4)(n-1)^{2}}{8}\Big]$;
\end{spacing}

\begin{spacing}{1.5}
\noindent $(3)$ $\displaystyle
 \Phi_{3}(\rho,n)=\frac{\rho^{n-6}}{8(n-2)}\Big\{8(n-2)\rho^{6}-(2n^{3}-11n^{2}+39n-44)\rho^{4}-2n^{2}\rho^{3}+\Big[\frac{13n^{2}}{4}+9(n-2)+
 \frac{17}{4}(n-5)(n-1)^{2}\Big]\rho^{2}+\frac{9}{4}n^{2}\rho-\frac{9}{4}(n-5)(n-1)^{2}\Big\}$.
\end{spacing}

For $\displaystyle \Phi_{1}(\rho,n)$, let $\displaystyle
g_{1}(\rho,n)=8(n-2)\rho^{4}-(2n^{3}-10n^{2}+34n-40)\rho^{2}+4(n-4)(n-1)^{2}$, then $\displaystyle
\Phi_{1}(\rho,n)=\frac{\rho^{n-4}}{8(n-2)}g_{1}(\rho,n)$.

Denote by $\rho_{1}\geq \rho_{2}\geq \rho_{3}\geq \rho_{4}$ the roots of $\displaystyle g_{1}(\rho,n)=0$. Since $\rho_{4}<0$, $\rho_{1}>0$ and $\displaystyle
g_{1}(0,n)=4(n-4)(n-1)^{2}>0$ for $n>4$. Then we have $\rho_{3}<0<\rho_{2}$. Moreover, as

\begin{center}
$\displaystyle g_{1}(\frac{n-1}{2},n)=\frac{n^{4}}{2}-3n^{3}-\frac{5n^{2}}{2}+12n-7>0$
\end{center}

\noindent for $n\geq7$ and

\begin{center}
$\displaystyle g_{1}(\frac{n-1.5}{2},n)=-\frac{3n^{3}}{8}-\frac{25n^{2}}{8}+\frac{93n}{32}+\frac{23}{16}<0$
\end{center}

\noindent for $n\geq2$, it follows that $\displaystyle \frac{n-1.5}{2}<\rho_{1}<\frac{n-1}{2}$, i.e., $\displaystyle \frac{n-1.5}{2}<\rho_{ag}(G_{1})<\frac{n-1}{2}$ for $n\geq7$.

The upper and lower bounds of $\displaystyle \rho_{ag}(G_{2})$ can also be obtained as above. Let
$\displaystyle g_{2}(\rho,n)=4(n-2)\rho^{4}-\frac{6n^{3}-31n^{2}+115n-136}{6}\rho^{2}-\frac{5n^{2}+5n}{6}\rho+\frac{2n^{2}}{3}
+\frac{19(n-4)(n-1)^{2}}{8}$, then $\displaystyle
\Phi_{2}(\rho,n)=\frac{\rho^{n-4}}{4(n-2)}g_{2}(\rho,n)$. Since
\begin{center}
$\displaystyle g_{2}(0,n)=\frac{2n^{2}}{3}+\frac{19(n-4)(n-1)^{2}}{8}>0$
\end{center}
\noindent for $n\geq4$.
Moreover,
\begin{center}
$\displaystyle g_{2}(\frac{n-1.5}{2},n)=\frac{n^{4}}{24}-\frac{43n^{3}}{48}-\frac{53n^{2}}{96}+\frac{143n}{64}+\frac{23}{32}>0$ \end{center}
for $n\geq22$, and
\begin{center}
$\displaystyle g_{2}(\frac{n-1.6}{2},n)=-\frac{n^{4}}{120}-\frac{17n^{3}}{30}-\frac{301n^{2}}{375}+\frac{22081n}{15000}+\frac{6487}{3750}<0$ \end{center}
for $n\geq2$.

It follows that $\displaystyle \frac{n-1.6}{2}<\rho_{ag}(G_{2})<\frac{n-1.5}{2}$ for $n\geq22$.

For $\displaystyle \Phi_{3}(\rho,n)$, let
\begin{center}
$\displaystyle
g_{3}(\rho,n)=8(n-2)\rho^{6}-(2n^{3}-11n^{2}+39n-44)\rho^{4}-2n^{2}\rho^{3}+\Big[\frac{13n^{2}}{4}+9(n-2)+
 \frac{17}{4}(n-5)(n-1)^{2}\Big]\rho^{2}+\frac{9}{4}n^{2}\rho-\frac{9}{4}(n-5)(n-1)^{2}$,
 \end{center}
 then $\displaystyle
\Phi_{3}(\rho,n)=\frac{\rho^{n-6}}{8(n-2)}g_{3}(\rho,n)$.
Clearly,
\begin{center}
$g_{3_{\rho}}^{(4)}(\rho,n)=(2880n-5760)\rho^{2}+24(-2n^{3}+11n^{2}-39n+44)$.
\end{center}
\noindent Since

\begin{center}
$\displaystyle g_{3_{\rho}}^{(4)}(\frac{n-1.6}{2},n)=672n^{3}-3480n^{2}+\frac{27576n}{5}-\frac{13152}{5}>0$
\end{center}

\noindent for $n\geq3$. Then $\displaystyle  g_{3_{\rho}}^{(4)}(\rho,n)>0$ for $\displaystyle \rho\geq\frac{n-1.6}{2}$ and $n\geq3$.

Thus if
$\displaystyle g_{3_{\rho}}^{(i)}(\frac{n-1.6}{2},n)>0$
for $0\leq i\leq 3$, then we can deduce that $\displaystyle
g_{3}(x,n)>0$ for $\displaystyle x\geq\frac{n-1.6}{2}$, i.e., $\displaystyle \rho_{ag}(G_{3})<\frac{n-1.6}{2}$.

We can verify that

$\displaystyle g_{3}(\frac{n-1.6}{2},n)=\frac{3n^{6}}{80}-\frac{149n^{5}}{200}+\frac{2293n^{4}}{1000}+\frac{4573n^{3}}{10000}-\frac{107861n^{2}}{50000}-\frac{54663n}{12500}-\frac{2619}{62500}>0$
for $n\geq17$.

$\displaystyle g_{3_{\rho}}^{(1)}(\frac{n-1.6}{2},n)=\frac{n^{6}}{2}-\frac{47n^{5}}{10}+\frac{1157n^{4}}{100}-\frac{1201n^{3}}{250}-\frac{1162n^{2}}{125}+\frac{79367n}{12500}+\frac{72954}{3125}>0$
 for $n\geq6$.

$\displaystyle g_{3_{\rho}}^{(2)}(\frac{n-1.6}{2},n)=9n^{5}-\frac{369n^{4}}{5}+\frac{9347n^{3}}{50}-\frac{3977n^{2}}{25}-\frac{5149n}{250}+\frac{15703}{250}>0$
for $n\geq5$.

$\displaystyle g_{3_{\rho}}^{(3)}(\frac{n-1.6}{2},n)=96n^{4}-\frac{3228n^{3}}{5}+\frac{6912n^{2}}{5}-\frac{26448n}{25}+\frac{3456}{25}>0$
for $n\geq4$.

Then it follows that $\displaystyle \rho_{ag}(G_{3})<\frac{n-1.6}{2}$ for $n\geq17$.
Since

$\displaystyle g_{3}(\frac{n-1.75}{2},n)=-\frac{47n^{5}}{128}+\frac{155n^{4}}{128}+\frac{1423n^{3}}{1024}-\frac{4115n^{2}}{2048}-\frac{43479n}{32768}-\frac{3105}{16384}<0$ for $n\geq4$.

Then (iii) follows. The proof is completed. \hfill$\Box$

\begin{lemma} \label{lemma 5}For $n\geq8$, $\displaystyle \frac{n-1}{2}<\rho_{ag}(S_{n}+e)<\frac{n}{2}$.
\end{lemma}
\noindent \textbf{Proof}  Since

\begin{equation*}       
  A_{ag}(S_{n}+e)=\begin{pmatrix}
    0 & \frac{n+1}{2\sqrt{2(n-1)}} & \frac{n+1}{2\sqrt{2(n-1)}} & \frac{n}{2\sqrt{n-1}} & \ldots &
    \frac{n}{2\sqrt{n-1}} &
    \frac{n}{2\sqrt{n-1}}\\  
    \frac{n+1}{2\sqrt{2(n-1)}} & 0 & 1 & 0 & \ldots & 0 & 0\\  
    \frac{n+1}{2\sqrt{2(n-1)}} & 1 & 0 & 0 & \ldots & 0 & 0\\
    \frac{n}{2\sqrt{n-1}} & 0 & 0 & 0 & \ldots & 0 & 0\\  
    \vdots & \vdots & \vdots & \vdots & \ddots & \vdots & \vdots\\
    \frac{n}{2\sqrt{n-1}} & 0 & 0 & 0 & \ldots & 0 & 0\\
    \frac{n}{2\sqrt{n-1}} & 0 & 0 & 0 & \ldots & 0 & 0\\

  \end{pmatrix}.                
\end{equation*}

\noindent Then its arithmetic-geometric characteristic polynomial
$\Phi(S_{n}+e,\rho)=\textrm{det}(\rho E_{n}-A_{ag}(S_{n}+e))$

\begin{equation*}       
=                
  \begin{vmatrix}
   \rho & -\frac{n+1}{2\sqrt{2(n-1)}} & -\frac{n+1}{2\sqrt{2(n-1)}} & -\frac{n}{2\sqrt{n-1}} & \ldots &
   -\frac{n}{2\sqrt{n-1}} &
   -\frac{n}{2\sqrt{n-1}}\\  
    -\frac{n+1}{2\sqrt{2(n-1)}} &\rho & -1 & 0 & \ldots & 0 & 0\\  
    -\frac{n+1}{2\sqrt{2(n-1)}} & -1 &\rho & 0 & \ldots & 0 & 0\\
    -\frac{n}{2\sqrt{n-1}} & 0 & 0 &\rho & \ldots & 0 & 0\\  
    \vdots & \vdots & \vdots & \vdots & \ddots & \vdots & \vdots\\
    -\frac{n}{2\sqrt{n-1}} & 0 & 0 & 0 & \ldots &\rho & 0\\
    -\frac{n}{2\sqrt{n-1}} & 0 & 0 & 0 & \ldots & 0 &\rho\\

  \end{vmatrix}              
\end{equation*}

$\displaystyle
=\rho^{n-4}(\rho+1)\Big[\rho^{3}-\rho^{2}-\frac{n^{3}-2n^{2}+2n+1}{4(n-1)}\rho+\frac{(n-3)n^{2}}{4(n-1)}\Big]$.

\noindent Let
\begin{center}
$\displaystyle
t_{1}(\rho,n)=\rho^{3}-\rho^{2}-\frac{n^{3}-2n^{2}+2n+1}{4(n-1)}\rho+\frac{(n-3)n^{2}}{4(n-1)}$.
\end{center}
\noindent Suppose that
$\displaystyle\rho_{1}\geq\rho_{2}\geq\rho_{3}$ are the three roots of $\displaystyle t_{1}(\rho,n)=0$. Obviously,
$\rho_{1}=\rho_{ag}(G)$. It remains to show that $\displaystyle
\frac{n-1}{2}<\rho_{1}<\frac{n}{2}$. Because for $n\geq8$,
\begin{center}
$\displaystyle t_{1}(0,n)=\frac{(n-3)n^{2}}{4(n-1)}>0$,
\end{center}

\begin{eqnarray*}
t_{1}\Big(\frac{n-1}{2},n\Big)&=&\frac{(n-1)^{4}-2(n-1)^{3}-(n-1)(n^{3}-2n^{2}+2n+1)+2n^{2}(n-3)}{8(n-1)}\\
&=&\frac{-n^{3}+2n^{2}-9n+4}{8(n-1)}<0,
\end{eqnarray*}
\noindent and

\begin{eqnarray*}
t_{1}\Big(\frac{n}{2},n\Big)&=&\frac{(n-1)n^{3}-2n^{2}(n-1)-n(n^{3}-2n^{2}+2n+1)+2n^{2}(n-3)}{8(n-1)}\\
&=&\frac{n^{3}-6n^{2}-n}{8(n-1)}>0.
\end{eqnarray*}

\noindent  It follows that $\displaystyle \frac{n-1}{2}<\rho_{1}<\frac{n}{2}$, i.e.,
$\displaystyle \frac{n-1}{2}<\rho_{ag}(S_{n}+e)<\frac{n}{2}$. \hfill$\Box$

\begin{theorem} Among all unicyclic graphs of order $n\geq5$, the cycle of order $n$ has the smallest $AG$ spectral radius and

\begin{itemize}
  \item[(i)] if $n=5,6,7$, then $S_{n}+e$, $G_{2}$, $G_{3}$ and $G_{1}$ are, respectively, the unique unicyclic graphs with the first four largest $AG$ spectral radii;
  \item[(ii)] if $8\leq n\leq 15$, then $S_{n}+e$, $G_{2}$, $G_{1}$ and $G_{3}$ are, respectively, the unique unicyclic graphs with the first four largest $AG$ spectral radii;
  \item[(iii)] if $16\leq n$, then $S_{n}+e$, $G_{1}$, $G_{2}$ and $G_{3}$ are, respectively, the unique unicyclic graphs with the first four largest $AG$ spectral radii.
\end{itemize}

\end{theorem}

\noindent \textbf{Proof}  For any unicyclic graph $G$ of order $n\geq5$.
Suppose that $\displaystyle h(x,y)=\frac{x+y}{2\sqrt{xy}}$, then
$\displaystyle h(x,y)\geq h(2,2)=1$ for $x,y\geq1$. By Lemmas
\ref{lemma 2} and \ref{lemma 3}, we
have $\displaystyle\rho_{ag}(G)\geq\rho(A(G)) \geq\rho(A(C_{n}))=\rho_{ag}(C_{n})=2.$ This leads to the cycle of order $n$ has the smallest $AG$ spectral radius.

\begin{figure}[H]
  \centering
    \includegraphics[width=9.5cm]{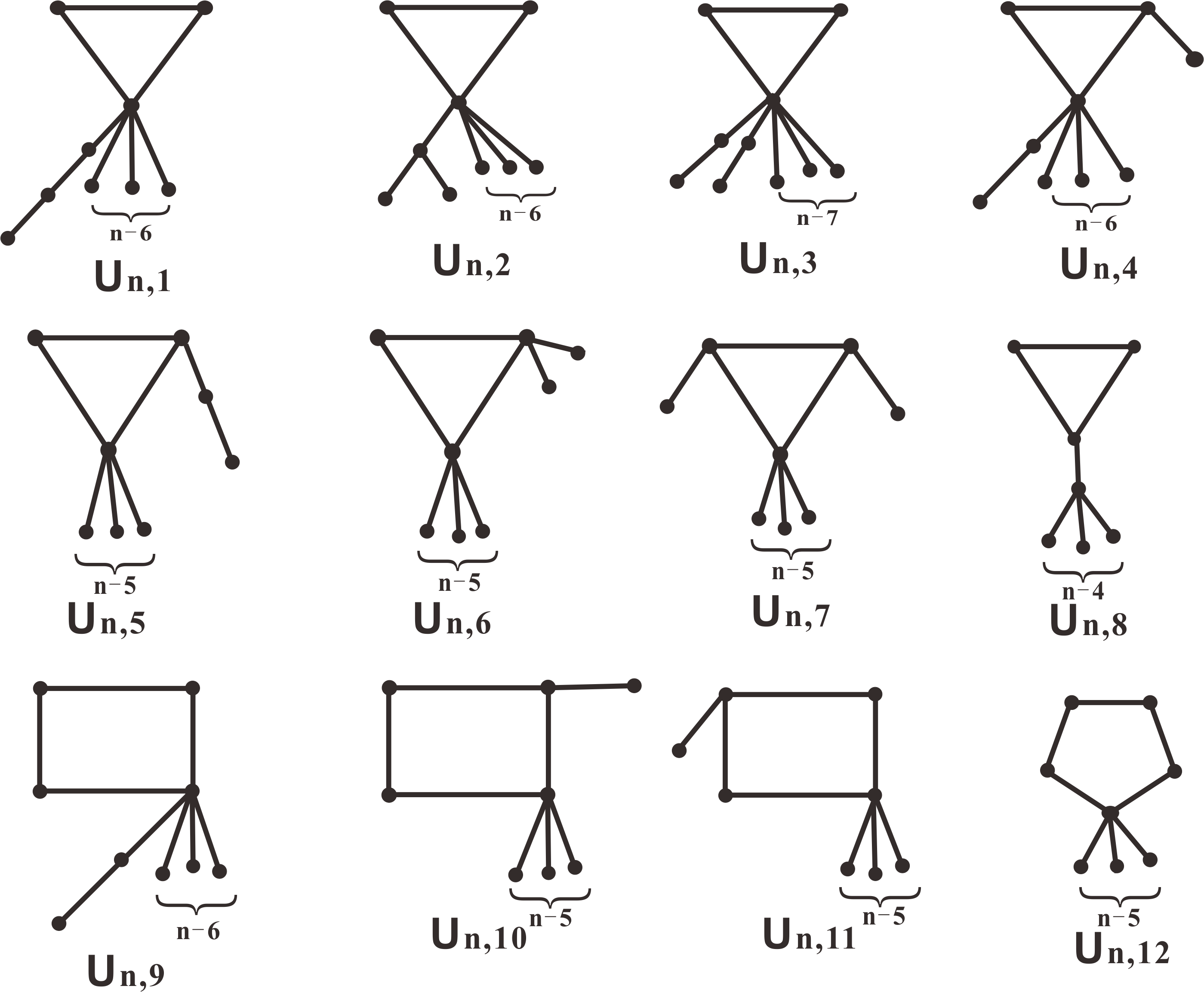}
  \caption{\small  The unicyclic graphs  ${U}_{n,i}$ with $i=1,...,12$.}
  \label{fig2}
\end{figure}

Denote by $\Delta$ the maximum degree of $G$, then $2\leq \Delta \leq n-1$. If $\Delta=n-1$, then $G\cong S_{n}+e$. If $\Delta=n-2$, then $G\cong G_{i}$ with $i=1,2,3$. If $\Delta=n-3$, then $G\cong U_{n,i}$ with $i=1,...,12$ are listed in Fig. \ref{fig2}.

For $5\leq n\leq7$, the approximate values of $\rho_{ag}$ of unicyclic graphs are given in Tables \ref{table 3}, \ref{table 4}, \ref{table 5} of the  Appendix 1. It is obvious that
$S_{n}+e$, $G_{2}$, $G_{3}$ and $G_{1}$ are, respectively, the unique unicyclic graphs with the first four largest $AG$ spectral radii.

For $8\leq n\leq15$, if $\Delta=n-3$, then the approximate values of $\rho_{ag}$ of unicyclic graphs are given in Table \ref{table 2} of the Appendix 1, and we found, in this case, $\displaystyle \rho_{ag}(G)<\frac{n-1.45}{2}$. Then by Lemmas \ref{lemma 7}, \ref{lemma 6} and \ref{lemma 5}, we obtain $S_{n}+e$, $G_{2}$, $G_{1}$ and $G_{3}$ are, respectively, the unique unicyclic graphs with the first four largest $AG$ spectral radii.

For $16\leq n\leq21$, if $G$ is not one of $S_{n}+e$, $G_{1}$, $G_{2}$ and $G_{3}$, then by Lemmas \ref{lemma 7}, \ref{lemma 6} and \ref{lemma 5}, we obtain

\begin{center}
$\displaystyle \rho_{ag}(G)<\frac{n-1.65}{2}<\rho_{ag}(G_{3})<\rho_{ag}(G_{1})<\rho_{ag}(G_{2})<\rho_{ag}(S_{n}+e)<\frac{n}{2}$.
\end{center}

\noindent Similarly, we have

\begin{center}
$\displaystyle \rho_{ag}(G)<\frac{n-1.75}{2}<\rho_{ag}(G_{3})<\rho_{ag}(G_{2})<\rho_{ag}(G_{1})<\rho_{ag}(S_{n}+e)<\frac{n}{2}$,
\end{center}

\noindent for $22\leq n $.

The results follow. \hfill$\Box$

\section{Conclusions}\label{sec5}
 \noindent

In this paper, we characterize the unicyclic graphs with the smallest and first four largest
arithmetic-geometric spectral radii. And further investigation of the following more general problems would be interesting.

\noindent \textbf{Problem 4.1} Determine more properties of the $AG$ matrix.

\noindent \textbf{Problem 4.2} What is the structure of the first two maximum bicyclic graphs for the $AG$ spectral radius? We conjecture the unique two bicyclic graphs are the bicyclic graphs of order $n$ with maximum degree $n-1$ for $n\geq7$.

\section{Declaration of competing interest}
\label{sec:3}

The authors declare that they have no known competing financial interests or personal
relationships that could have appeared to influence the work reported in this paper.

\section{Acknowledgments}
\label{sec:5}

This work is supported by NSFC (No. 12171402) and the Fundamental Research Funds for the Central Universities (No. 20720190062).

\bibliographystyle{abbrv}
\bibliography{mybibfile}

\newpage

\section*{Appendix 1}
\label{sec:Appendix 1}

\begin{table}[H]
\centering
\caption{The approximate value of $\rho_{ag}$ of $G_{i}$ with $i=1,2,3$, $ \frac{n-1.65}{2}$, $\frac{n-1}{2}$ for $5\leq n\leq 21$. }
\renewcommand\arraystretch{1}
\setlength\tabcolsep{14pt}
  \begin{tabular}{  c | l  l l l  l l l }
    \hline
     $n$ & $\displaystyle\rho_{ag}(G_{1})$ & $\displaystyle\rho_{ag}(G_{2})$ & $\displaystyle\rho_{ag}(G_{3})$& $ \frac{n-1.65}{2}$ & $\frac{n-1.45}{2}$ & $ \frac{n-1}{2}$\\ \hline
     8 & 3.3765 & 3.4571 & 3.3755& 3.175 & 3.275 & 3.5 \\ \hline
     9 & 3.8437 & 3.9009 & 3.8245& 3.675  & 3.775 &4 \\ \hline
     10 & 4.3229 & 4.3632 & 4.2890 & 4.175&  4.275 & 4.5  \\ \hline
     11 & 4.8087 & 4.8369 & 4.7633& 4.675 & 4.775 &5  \\ \hline
     12 & 5.2987 & 5.3176 & 5.2440& 5.175 & 5.275 & 5.5 \\ \hline
     13 & 5.7913 & 5.8031 & 5.7292& 5.675 & 5.775 & 6  \\ \hline
     14 & 6.2856 & 6.2917 & 6.2174& 6.175 & 6.275 & 6.5  \\ \hline
     15 & 6.7812 & 6.7827& 6.7079 & 6.675 & 6.775 &7 \\ \hline
     16 & 7.2777 & 7.2753& 7.2001  & 7.175  & 7.275 &7.5 \\ \hline
     17 & 7.7748 & 7.7693& 7.6936 & 7.675  & 7.775 &8 \\ \hline
     18 & 8.2725 & 8.2642& 8.1881 & 8.175 & 8.275 &8.5 \\ \hline
     19 & 8.7705 & 8.7598& 8.6833 & 8.675  & 8.775 &9 \\ \hline
     20 & 9.2688 & 9.2561& 9.1792 & 9.175  & 9.275 &9.5 \\ \hline
     21 & 9.7673 & 9.7528& 9.6757 & 9.675  & 9.775 &10 \\
    \hline
     \end{tabular}
\label{table 1}
\end{table}

\begin{table}[H]
\caption{The approximate value of $\rho_{ag}$ of $U_{n,i}$ with $1\leq i\leq12$ for $8\leq n\leq 15$}
\renewcommand\arraystretch{1}
  \begin{tabular}{  c | l l l l l l l l}
    \hline
      & 8 & 9 & 10 & 11 & 12 & 13 & 14 & 15       \\ \hline
     $\displaystyle\rho_{ag}(U_{n,1})$ & 2.9523 & 3.3751 & 3.8239 & 4.2884 & 4.7627 & 5.2435 & 5.7287& 6.217\\ \hline
     $\displaystyle\rho_{ag}(U_{n,2})$ & 2.9702 & 3.374 & 3.8109 & 4.2678 & 4.7371 & 5.2145 & 5.6972 &  6.1837\\ \hline
     $\displaystyle\rho_{ag}(U_{n,3})$ & 2.8938 & 3.2991 & 3.7355 & 4.1914 & 4.6597 & 5.1362 & 5.6182 & 6.1041 \\ \hline
     $\displaystyle\rho_{ag}(U_{n,4})$ & 2.9894 & 3.3827 & 3.8135 & 4.2671 & 4.7347 & 5.2111 & 5.6933& 6.1796 \\ \hline
     $\displaystyle\rho_{ag}(U_{n,5})$ & 3.0375 & 3.451 & 3.8963 & 4.3599 & 4.8343 & 5.3156 & 5.8015& 6.2904 \\ \hline
     $\displaystyle\rho_{ag}(U_{n,6})$ & 3.1534 & 3.5152 & 3.9296 & 4.375 & 4.8384 & 5.3127 & 5.7939 & 6.2797 \\ \hline
     $\displaystyle\rho_{ag}(U_{n,7})$ & 3.0673 & 3.4533 & 3.882 & 4.3355 & 4.8039 & 5.2813 & 5.7646 & 6.2519 \\ \hline
     $\displaystyle\rho_{ag}(U_{n,8})$ & 2.9753 & 3.4221 & 3.8954 & 4.3801 & 4.8704 & 5.3639 & 5.8591 & 6.3556 \\ \hline
     $\displaystyle\rho_{ag}(U_{n,9})$ & 2.8726 & 3.2984 & 3.7532 & 4.2241 & 4.7045 & 5.1906 & 5.6803 & 6.1726 \\ \hline
     $\displaystyle\rho_{ag}(U_{n,10})$ & 2.9553 & 3.3748 & 3.8283 & 4.2995 & 4.7806 & 5.2674 & 5.7579 & 6.2508 \\ \hline
     $\displaystyle\rho_{ag}(U_{n,11})$ & 2.9662 & 3.394 & 3.8534 & 4.3287 & 4.8125 & 5.3013 & 5.7932 & 6.2871 \\ \hline
     $\displaystyle\rho_{ag}(U_{n,12})$ & 2.8974 & 3.3498 & 3.8232 & 4.3068 & 4.7959 & 5.2883 & 5.7827& 6.2784 \\ \hline
     $\displaystyle\frac{n-1.45}{2}$ & 3.275 & 3.775 & 4.275 & 4.775 & 5.275 & 5.775 & 6.275& 6.775 \\
    \hline
     \end{tabular}
\label{table 2}
\end{table}

  \begin{table}[H]
   \caption{The approximate value of $\rho_{ag}$ of the unicyclic graph with order $n=5$}
\renewcommand\arraystretch{2}
  \centering
  \begin{tabular}{ | c | l | l |l | }
    \hline
The unicyclic graph G & $\displaystyle\rho_{ag}(G)$ & The unicyclic graph G & $\displaystyle\rho_{ag}(G)$\\
    \hline
    \begin{minipage}[b]{0.3\columnwidth}
		\centering
		\raisebox{-.35\height}{\includegraphics[width=1.4cm]{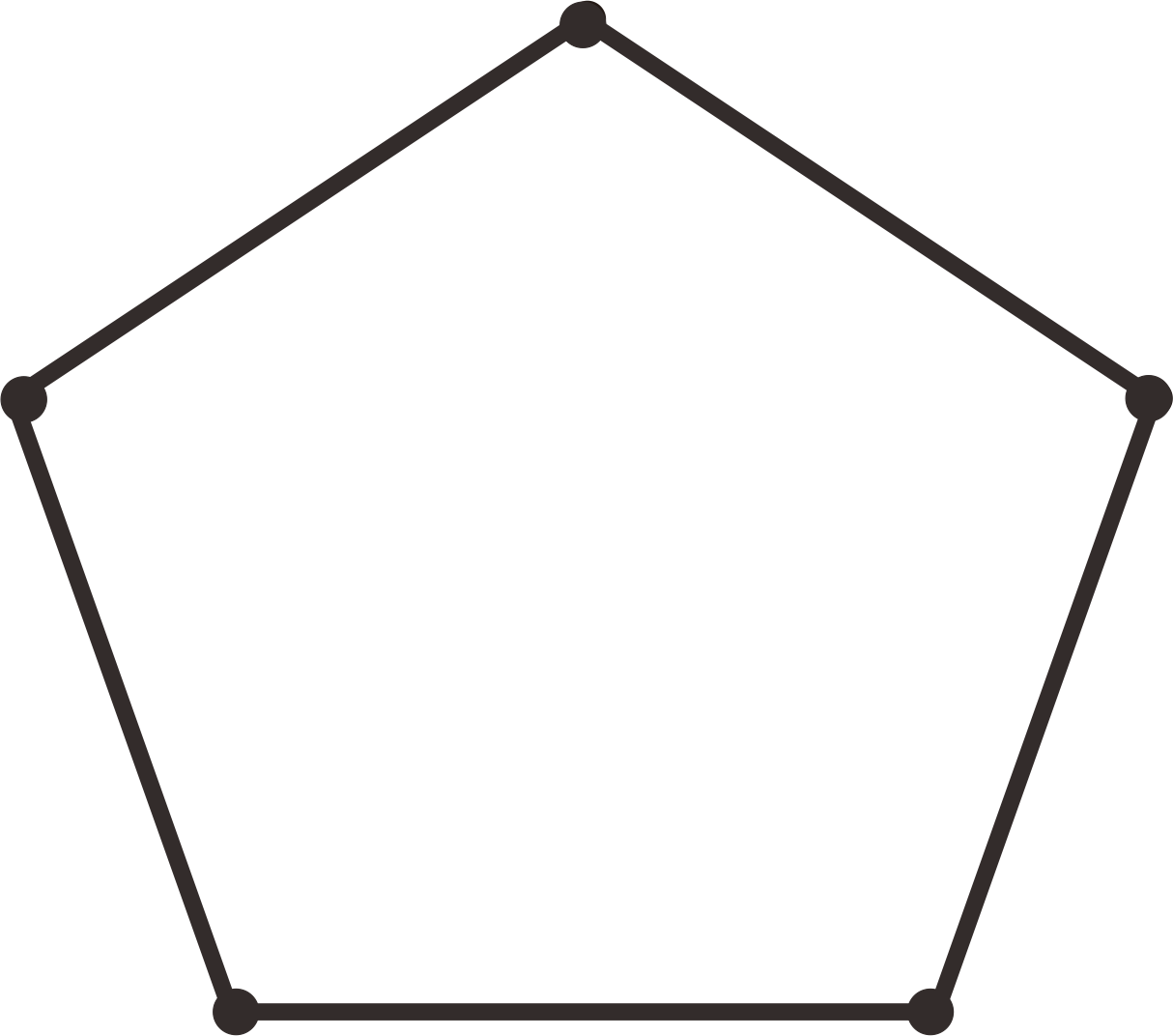}}
	\end{minipage}
    & \ \ \ \ \ 2
    & \begin{minipage}[b]{0.3\columnwidth}
		\centering
		\raisebox{-.35\height}{\includegraphics[width=1.75cm]{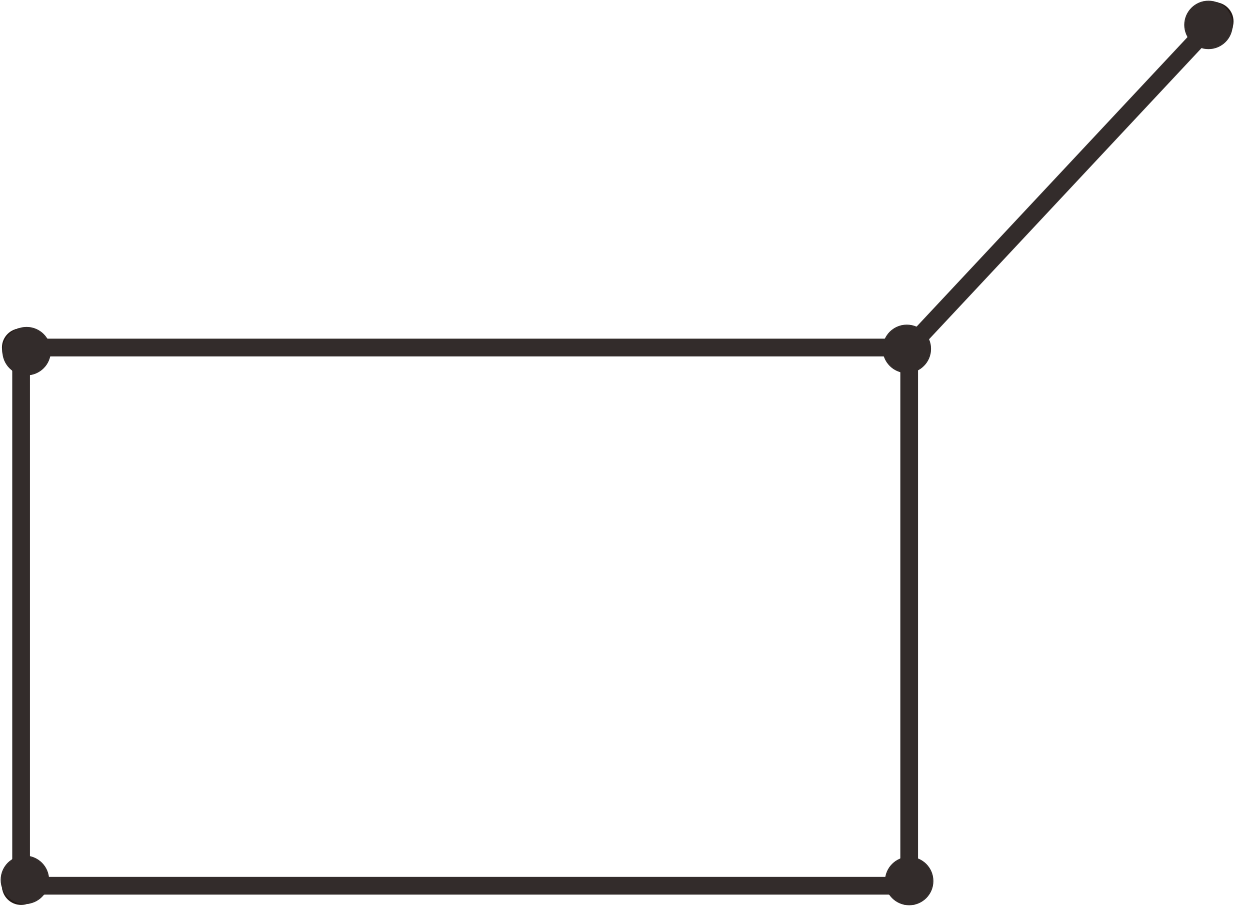}}
	\end{minipage}
    &  2.2066
    \\[3pt]
    \hline
   \begin{minipage}[b]{0.3\columnwidth}
		\centering
		\raisebox{-.30\height}{\includegraphics[width=2.3cm]{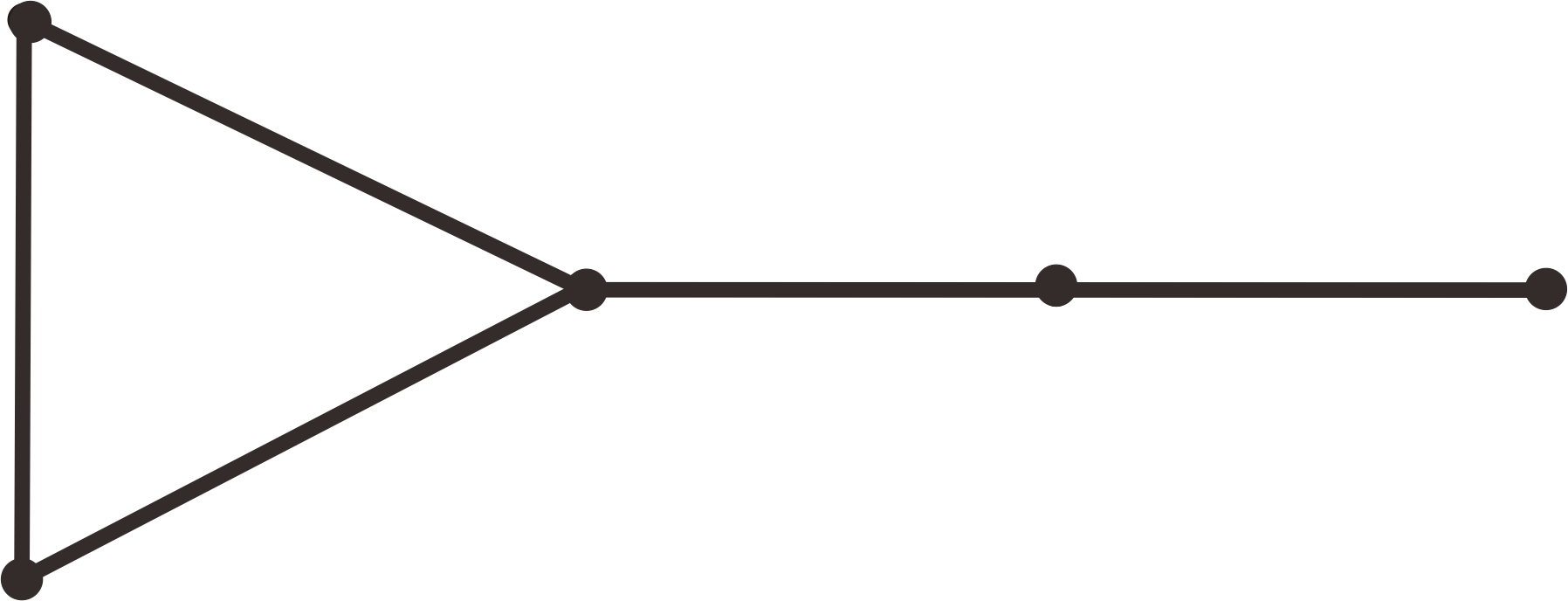}}
	\end{minipage}
    & 2.2543
    & \begin{minipage}[b]{0.3\columnwidth}
		\centering
		\raisebox{-.35\height}{\includegraphics[width=1.5cm]{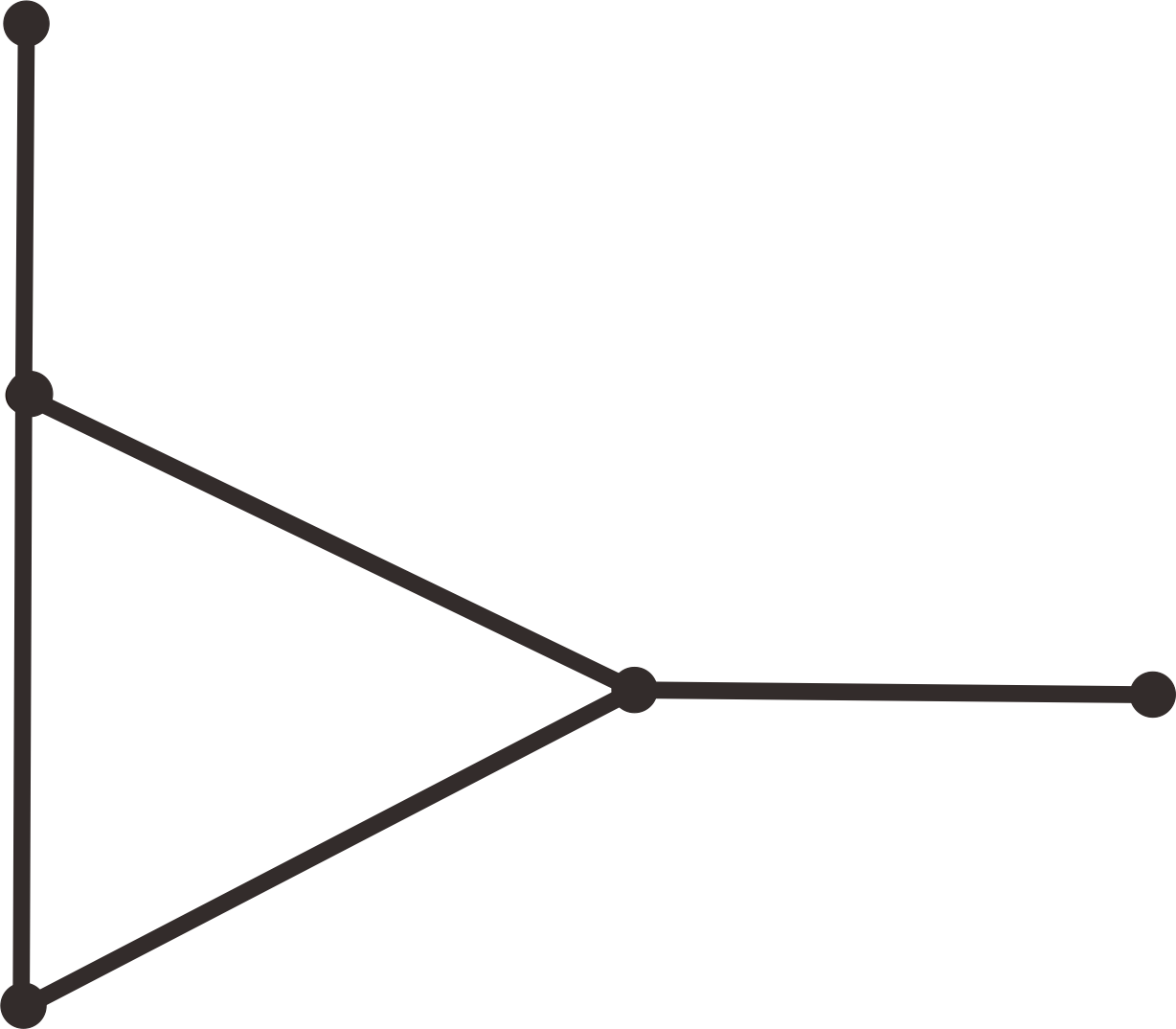}}
	\end{minipage}
    &  2.4149
    \\[3pt]
    \hline
  \begin{minipage}[b]{0.3\columnwidth}
		\centering
		\raisebox{-.30\height}{\includegraphics[width=1.7cm]{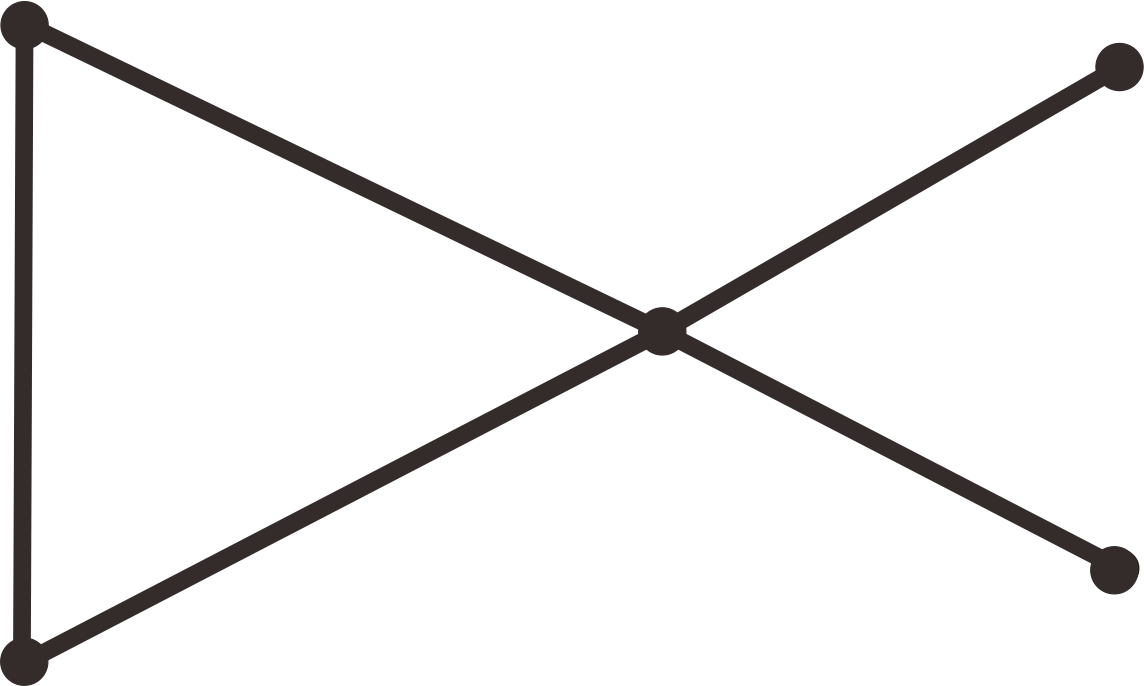}}
	\end{minipage}
    & 2.6035
    &
    &
    \\ \hline
\end{tabular}
\label{table 3}
\end{table}

\begin{table}[H]
  \caption{The approximate value of $\rho_{ag}$ of the unicyclic graph with order $n=6$}
\renewcommand\arraystretch{2}
  \centering
  \begin{tabular}{ | c | l | l |l |}
    \hline
   The unicyclic graph G & $\displaystyle\rho_{ag}(G)$ & The unicyclic graph G & $\displaystyle\rho_{ag}(G)$\\
    \hline
    \begin{minipage}[b]{0.3\columnwidth}
		\centering
		\raisebox{-.25\height}{\includegraphics[width=1cm]{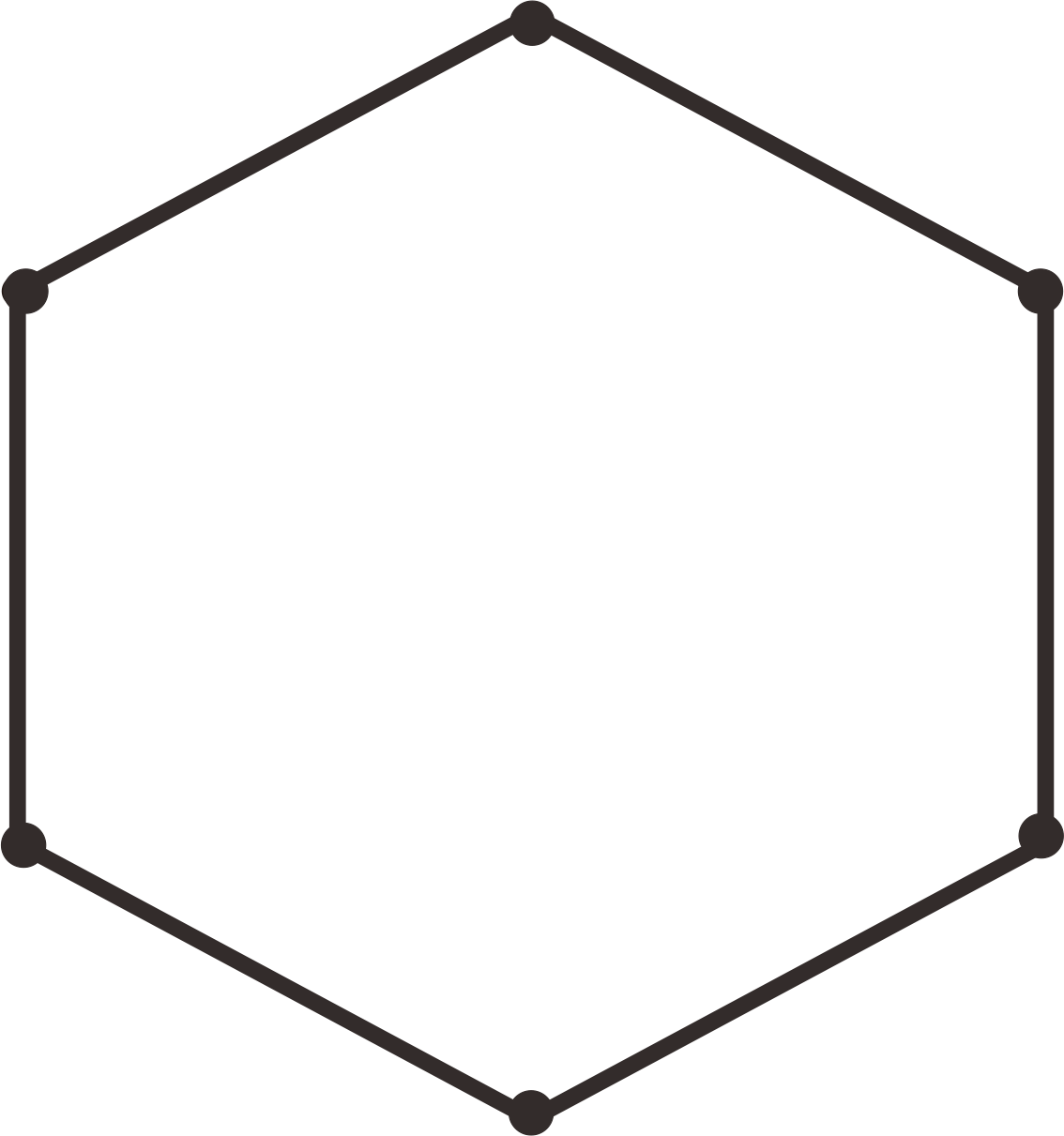}}
	\end{minipage}
    & \ \ \ \ 2
    & \begin{minipage}[b]{0.3\columnwidth}
		\centering
		\raisebox{-.28\height}{\includegraphics[width=1cm]{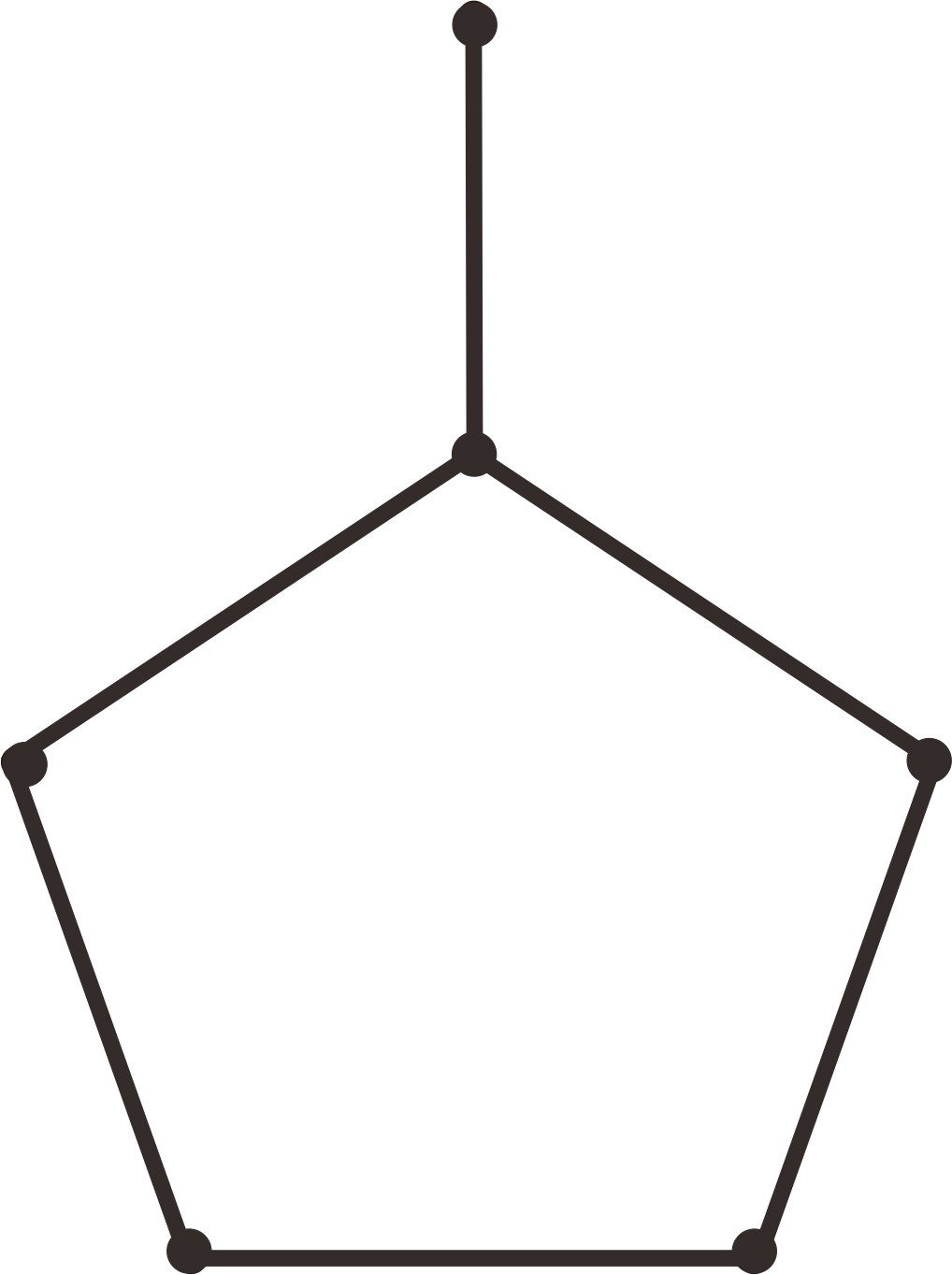}}
	\end{minipage}
    &  2.1785
    \\[3pt]
    \hline
    \begin{minipage}[b]{0.3\columnwidth}
		\centering
		\raisebox{-.28\height}{\includegraphics[width=1.5cm]{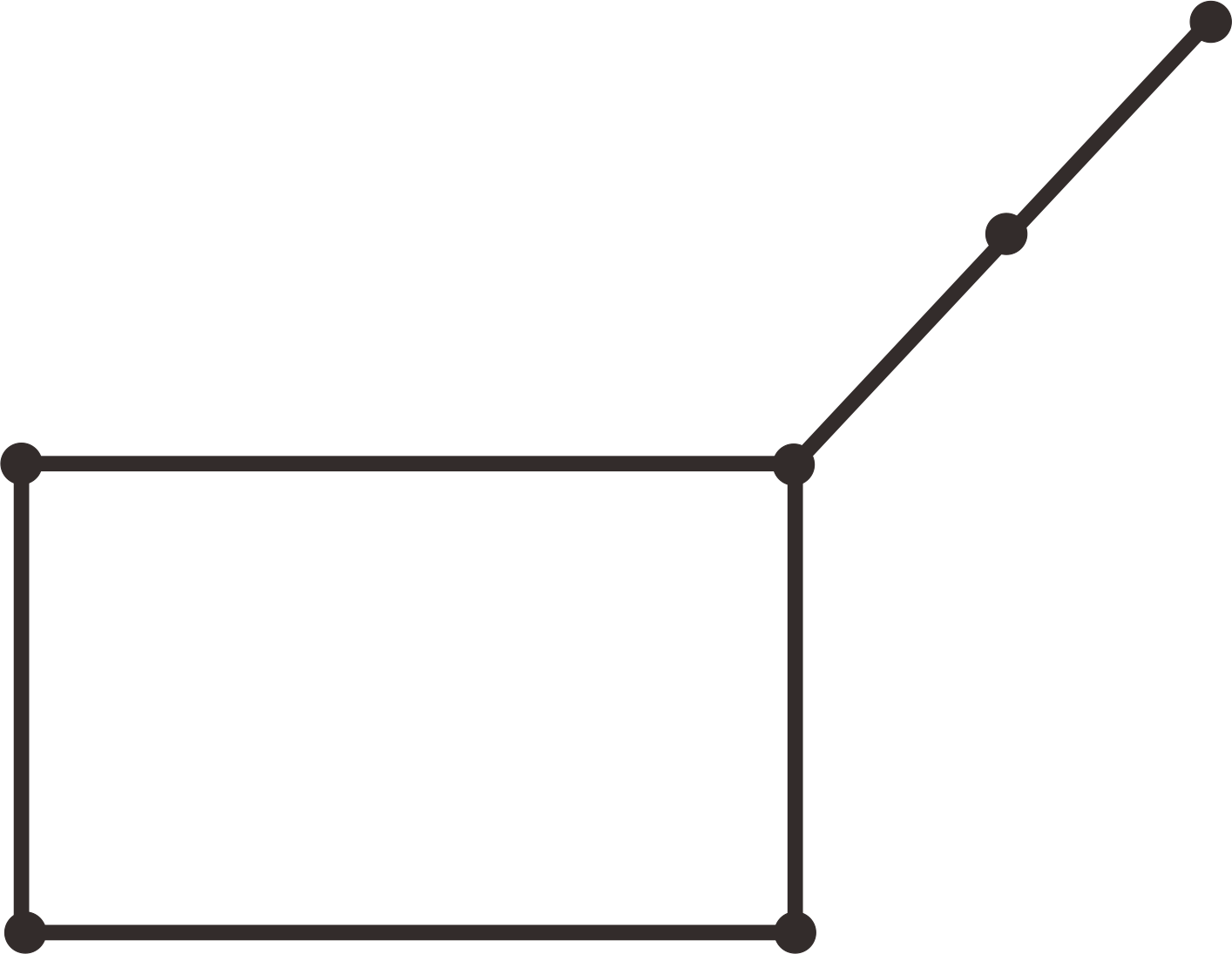}}
	\end{minipage}
    & 2.2096
    & \begin{minipage}[b]{0.3\columnwidth}
		\centering
		\raisebox{-.28\height}{\includegraphics[width=2cm]{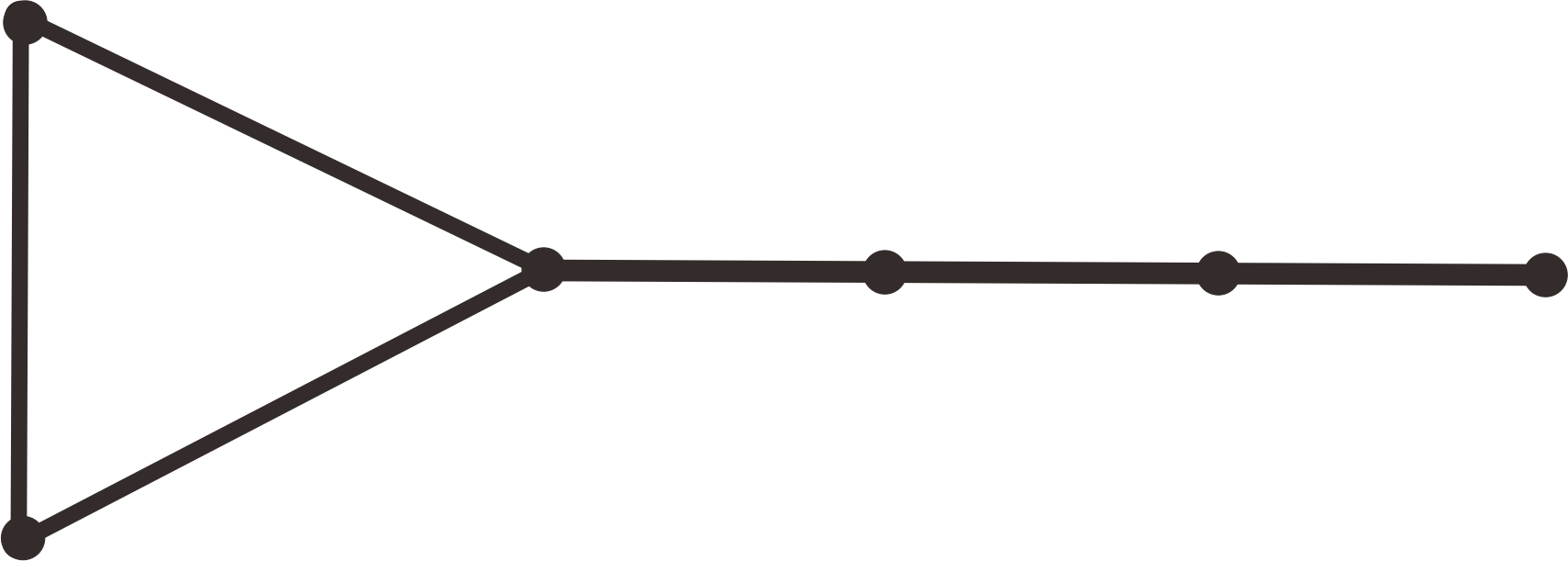}}
	\end{minipage}
    &  2.2632
    \\[3pt]
    \hline
 \begin{minipage}[b]{0.3\columnwidth}
		\centering
		\raisebox{-.33\height}{\includegraphics[width=1.5cm]{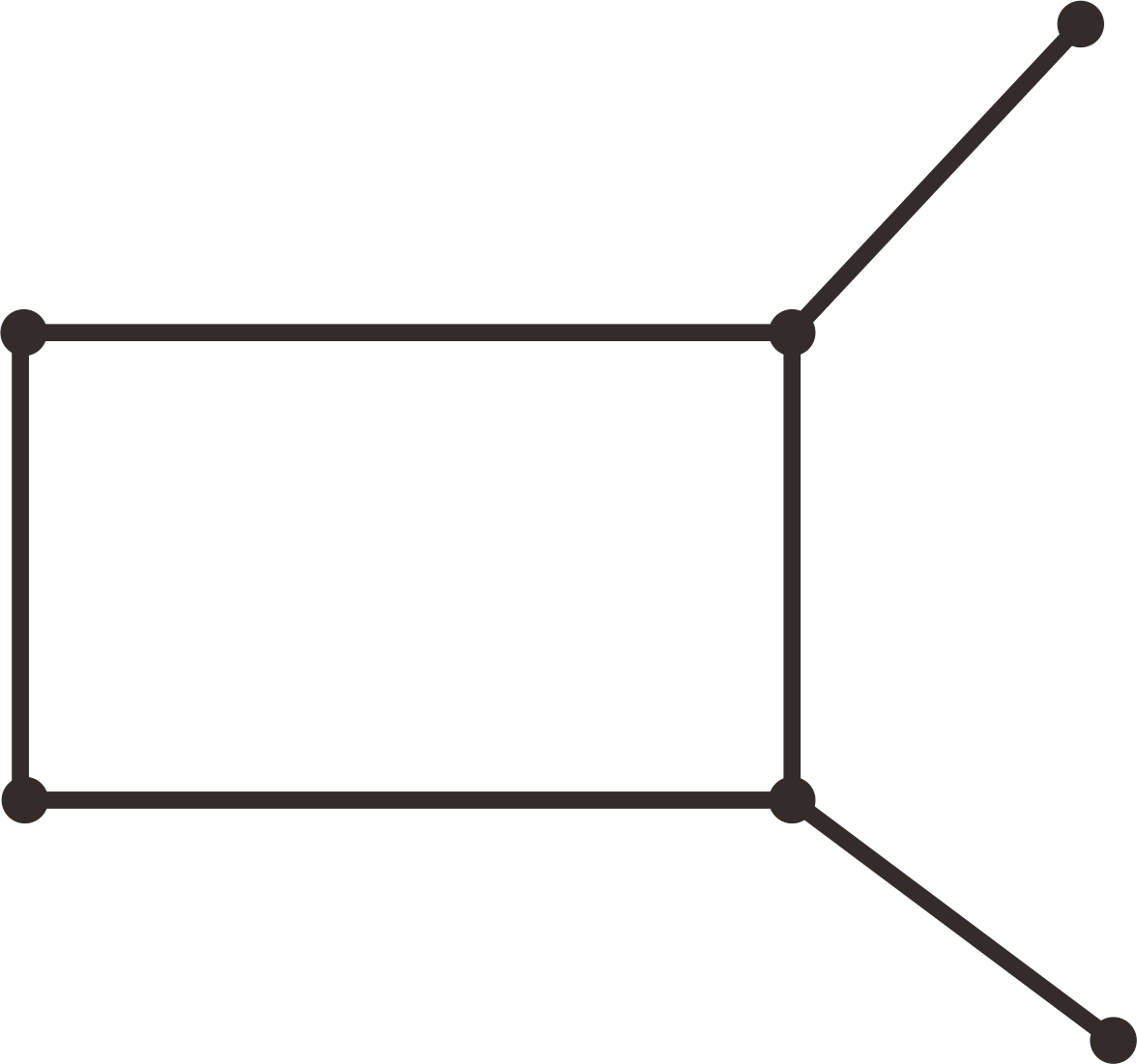}}
	\end{minipage}
    & 2.3439
    & \begin{minipage}[b]{0.3\columnwidth}
		\centering
		\raisebox{-.30\height}{\includegraphics[width=1.6cm]{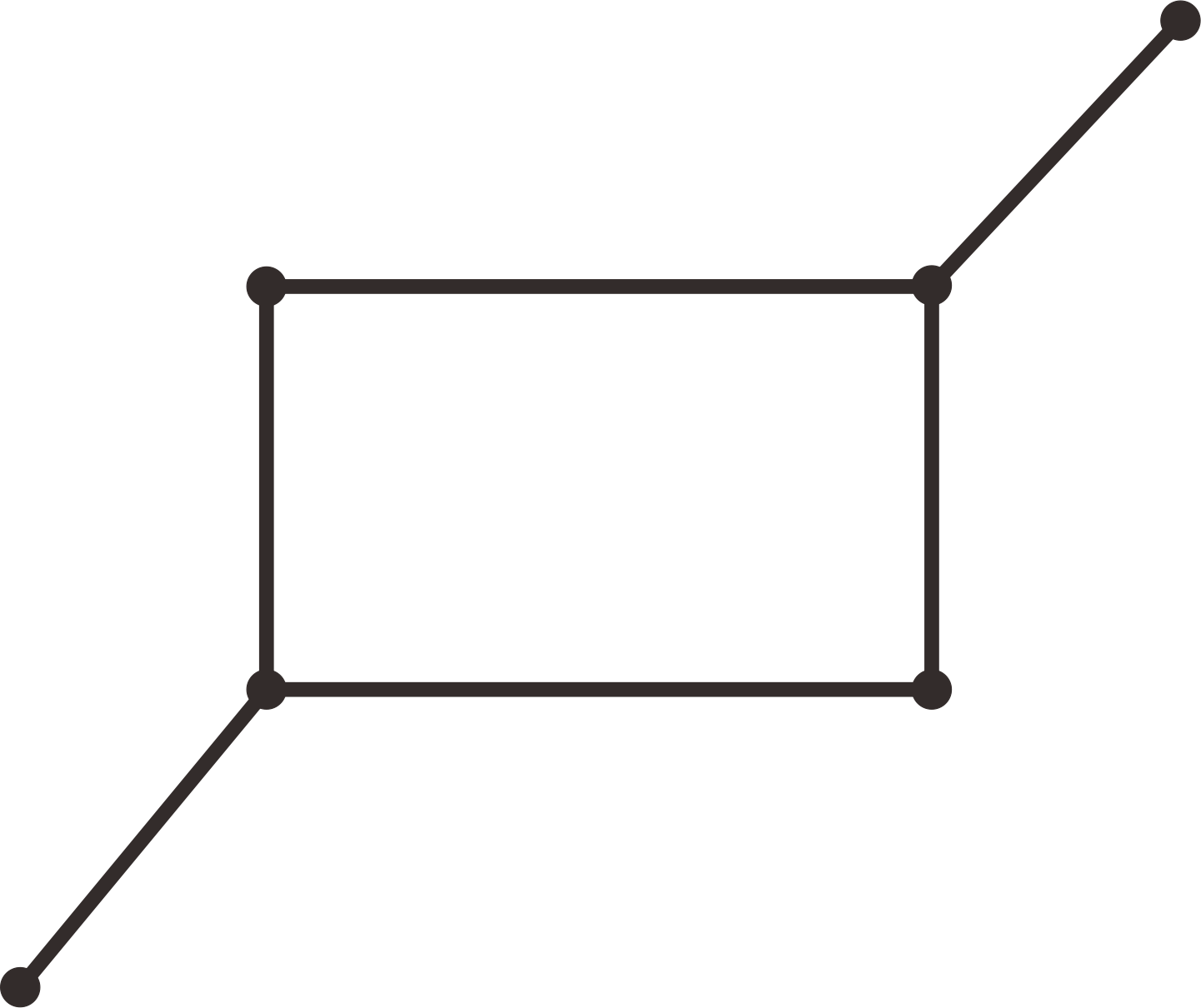}}
	\end{minipage}
    &  2.3452
    \\[3pt]
    \hline
     \begin{minipage}[b]{0.3\columnwidth}
		\centering
		\raisebox{-.25\height}{\includegraphics[width=1.6cm]{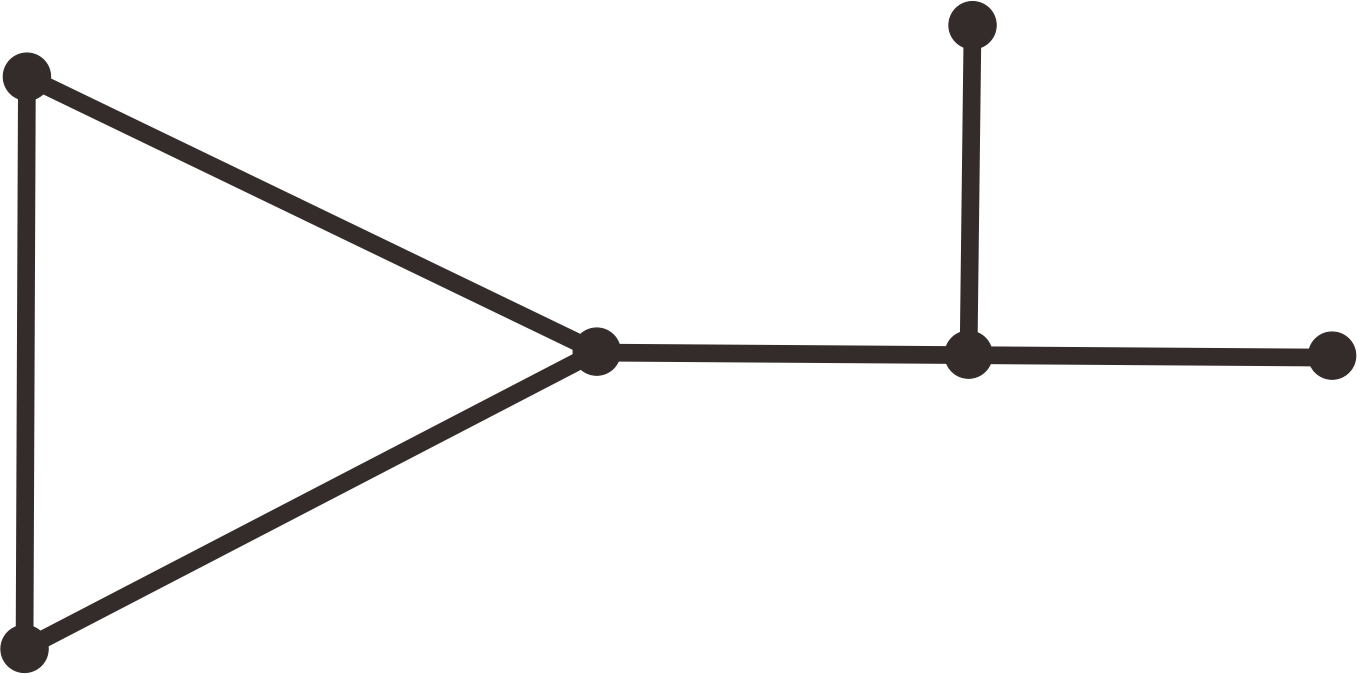}}
	\end{minipage}
    & 2.3551
    & \begin{minipage}[b]{0.3\columnwidth}
		\centering
		\raisebox{-.25\height}{\includegraphics[width=1.5cm]{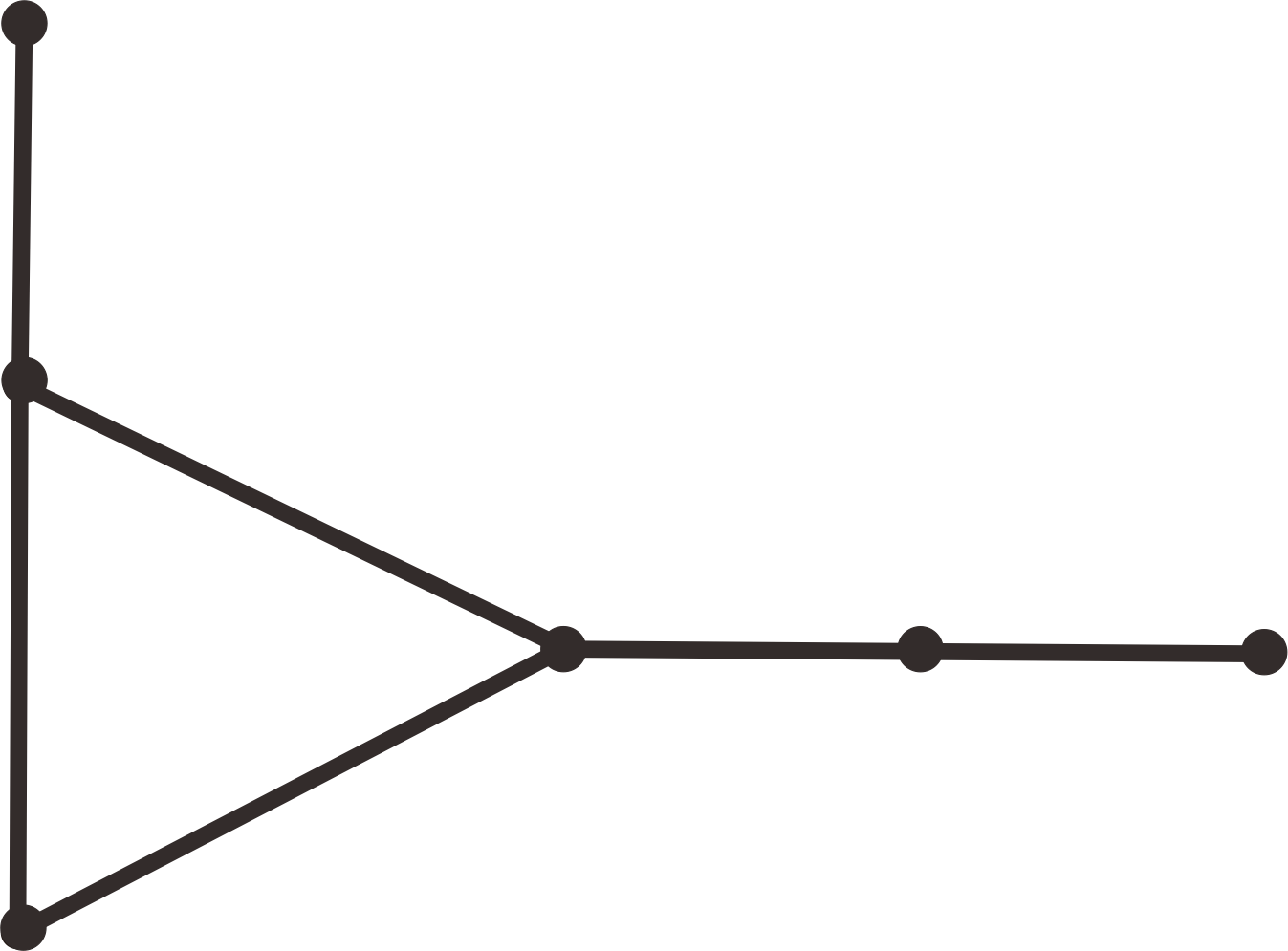}}
	\end{minipage}
    &  2.4095
    \\[3pt]
    \hline
     \begin{minipage}[b]{0.3\columnwidth}
		\centering
		\raisebox{-.35\height}{\includegraphics[width=1cm]{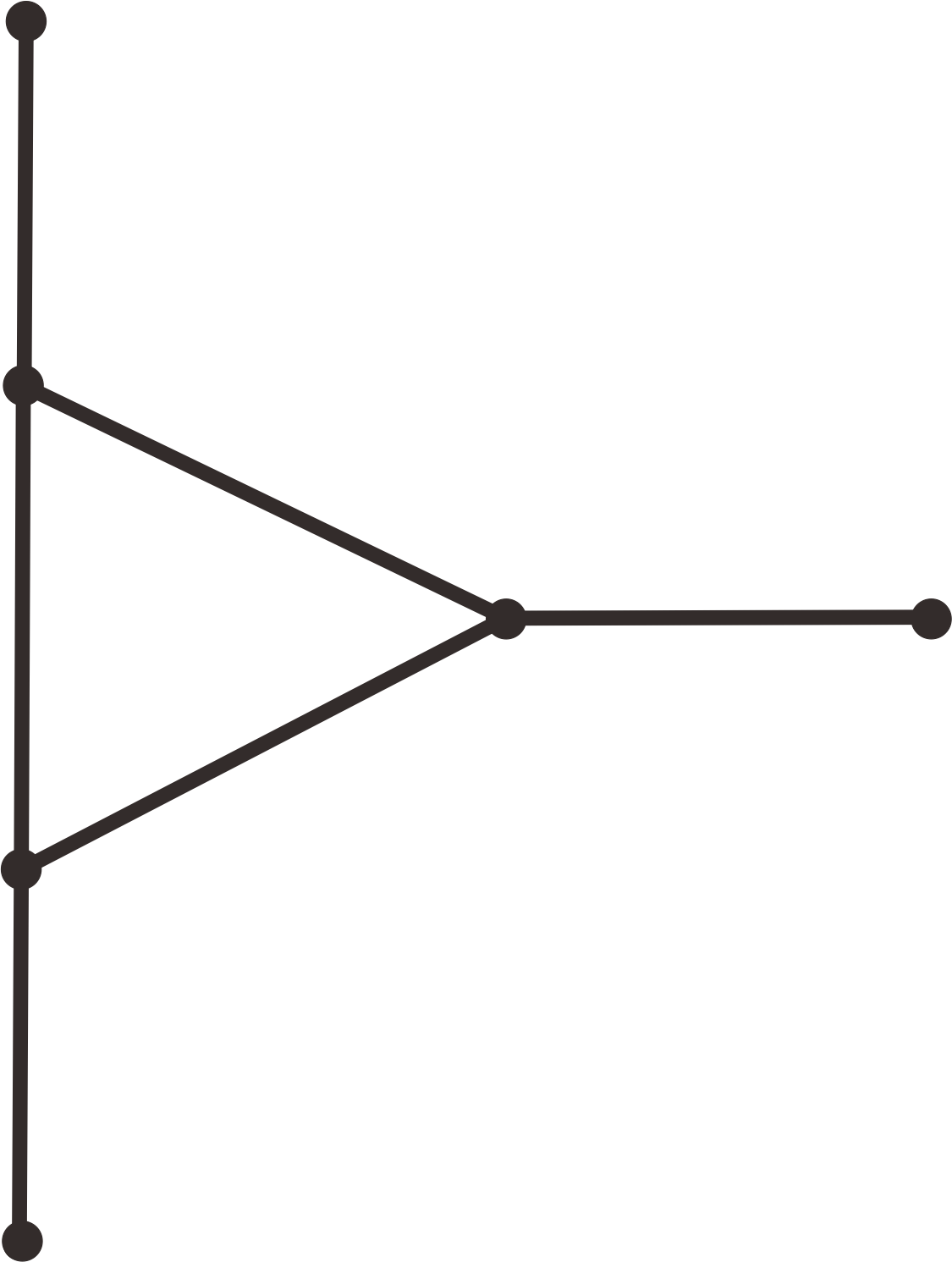}}
	\end{minipage}
    & 2.5275
    & \begin{minipage}[b]{0.3\columnwidth}
		\centering
		\raisebox{-.25\height}{\includegraphics[width=1.5cm]{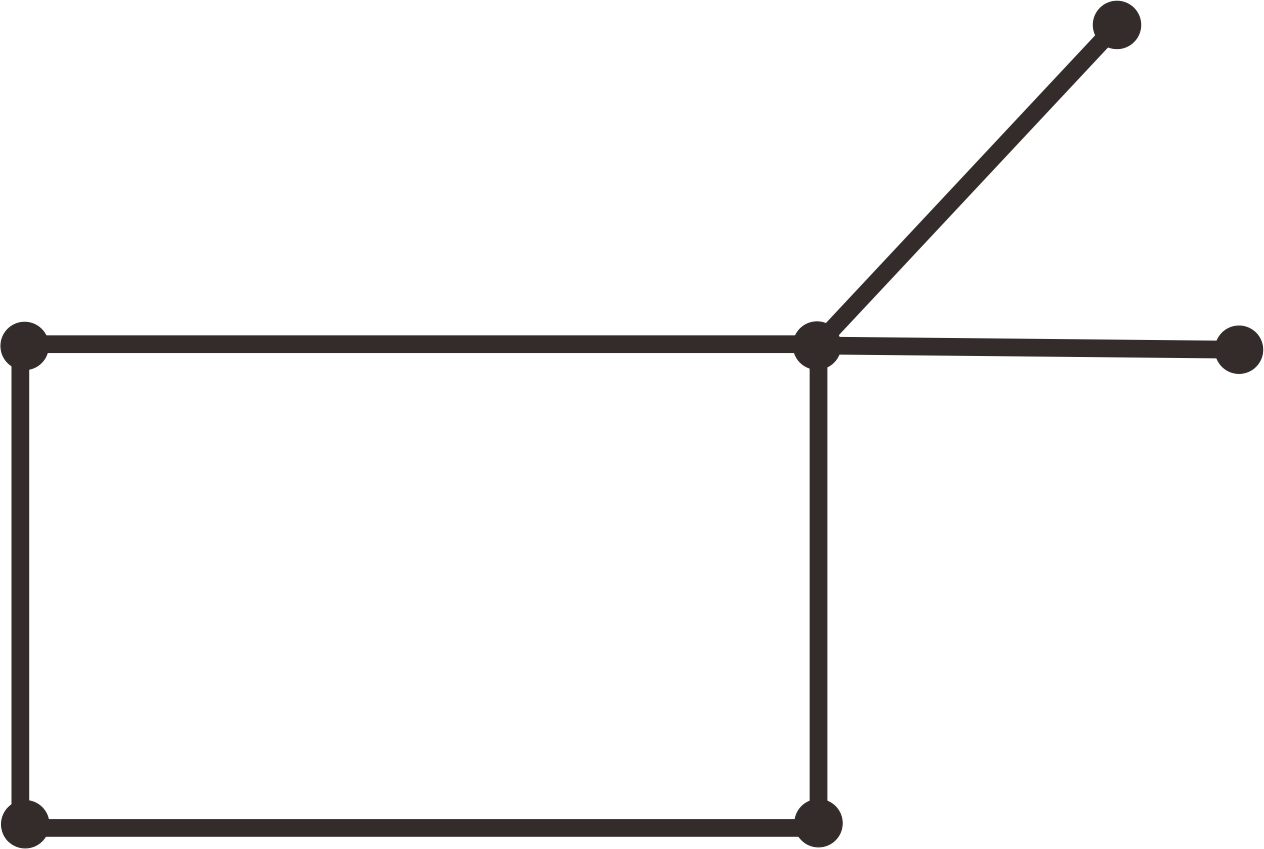}}
	\end{minipage}
    &  2.5295
    \\[3pt]
    \hline
     \begin{minipage}[b]{0.3\columnwidth}
		\centering
		\raisebox{-.25\height}{\includegraphics[width=1.6cm]{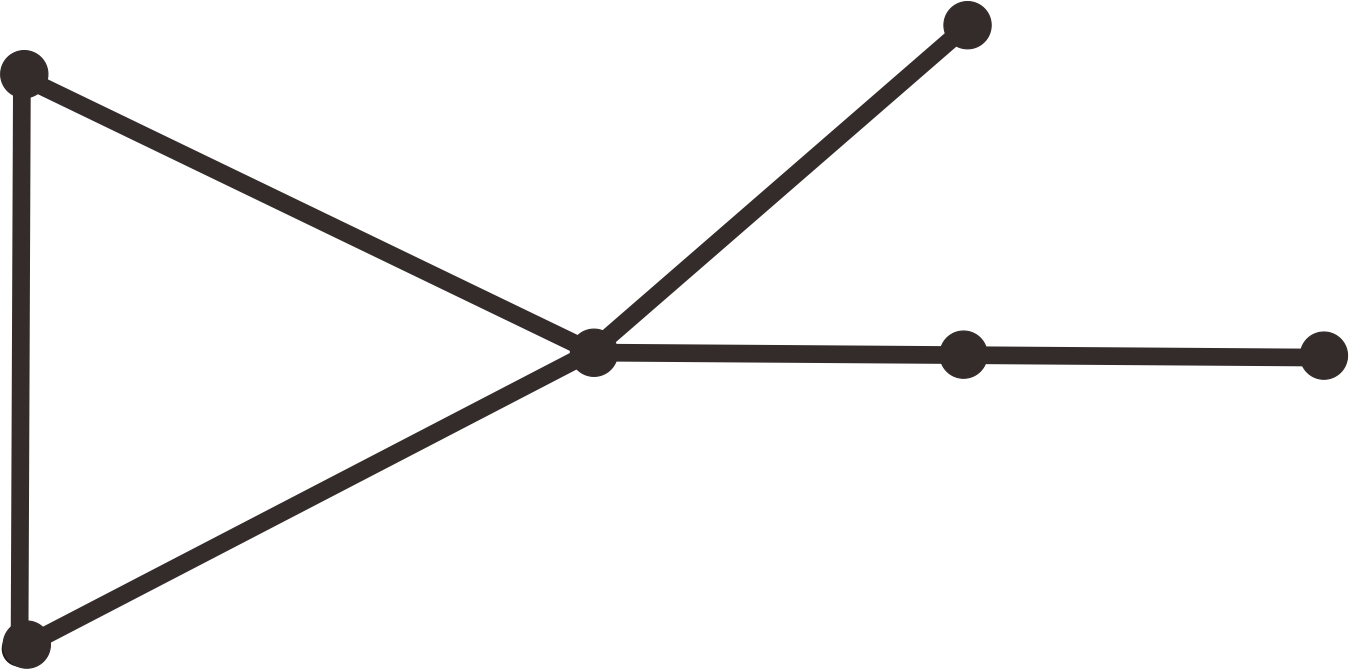}}
	\end{minipage}
    & 2.5695
    & \begin{minipage}[b]{0.3\columnwidth}
		\centering
		\raisebox{-.33\height}{\includegraphics[width=1.3cm]{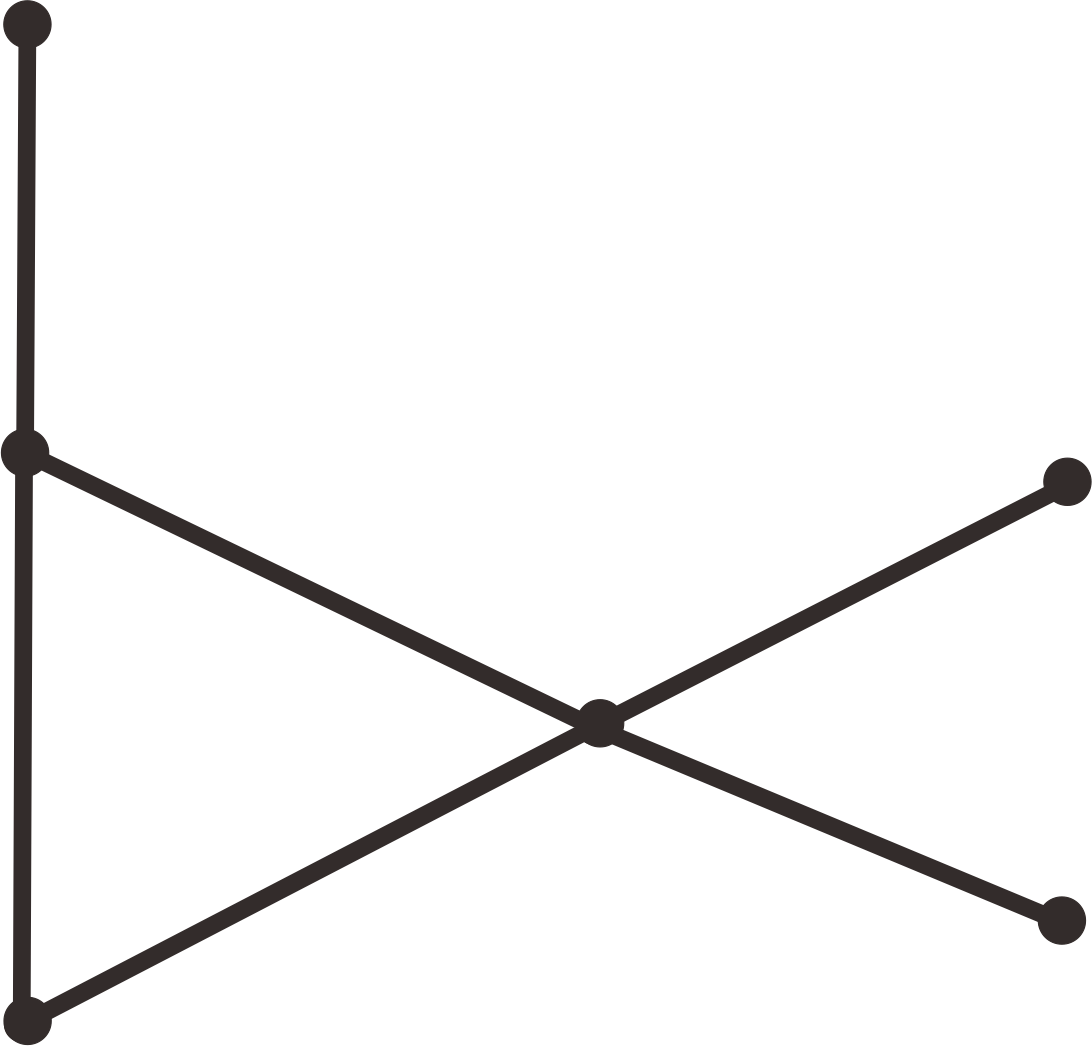}}
	\end{minipage}
    &  2.6879
    \\[3pt]
    \hline
     \begin{minipage}[b]{0.3\columnwidth}
		\centering
		\raisebox{-.30\height}{\includegraphics[width=1.3cm]{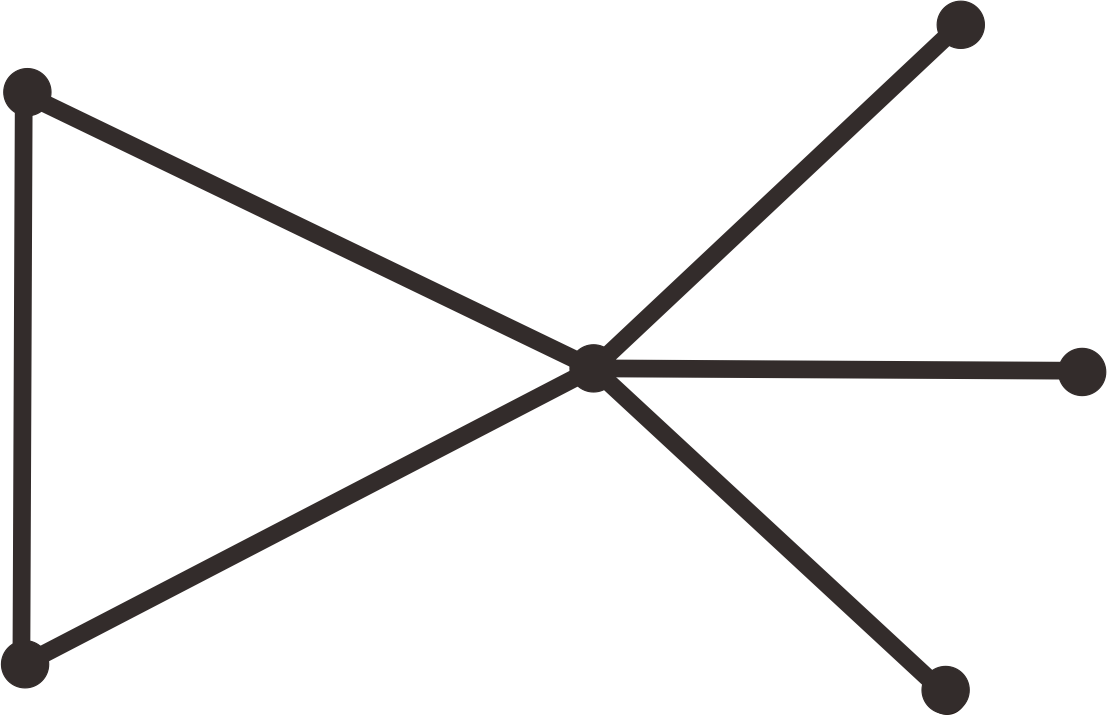}}
	\end{minipage}
    &  3.0113
    &
    &
     \\[3pt]
 \hline
\end{tabular}
\label{table 4}
\end{table}

\begin{table}[H]
  \caption{The approximate value of $\rho_{ag}$ of the unicyclic graph with order $n=7$}
\scalebox{0.1}{}
\renewcommand\arraystretch{1.3}
  \centering
  \begin{tabular}{ | c | l | l |l |}
    \hline
   The unicyclic graph G & $\displaystyle\rho_{ag}(G)$ & The unicyclic graph G & $\displaystyle\rho_{ag}(G)$\\
    \hline
     \begin{minipage}[b]{0.3\columnwidth}
		\centering
		\raisebox{-.25\height}{\includegraphics[width=0.9cm]{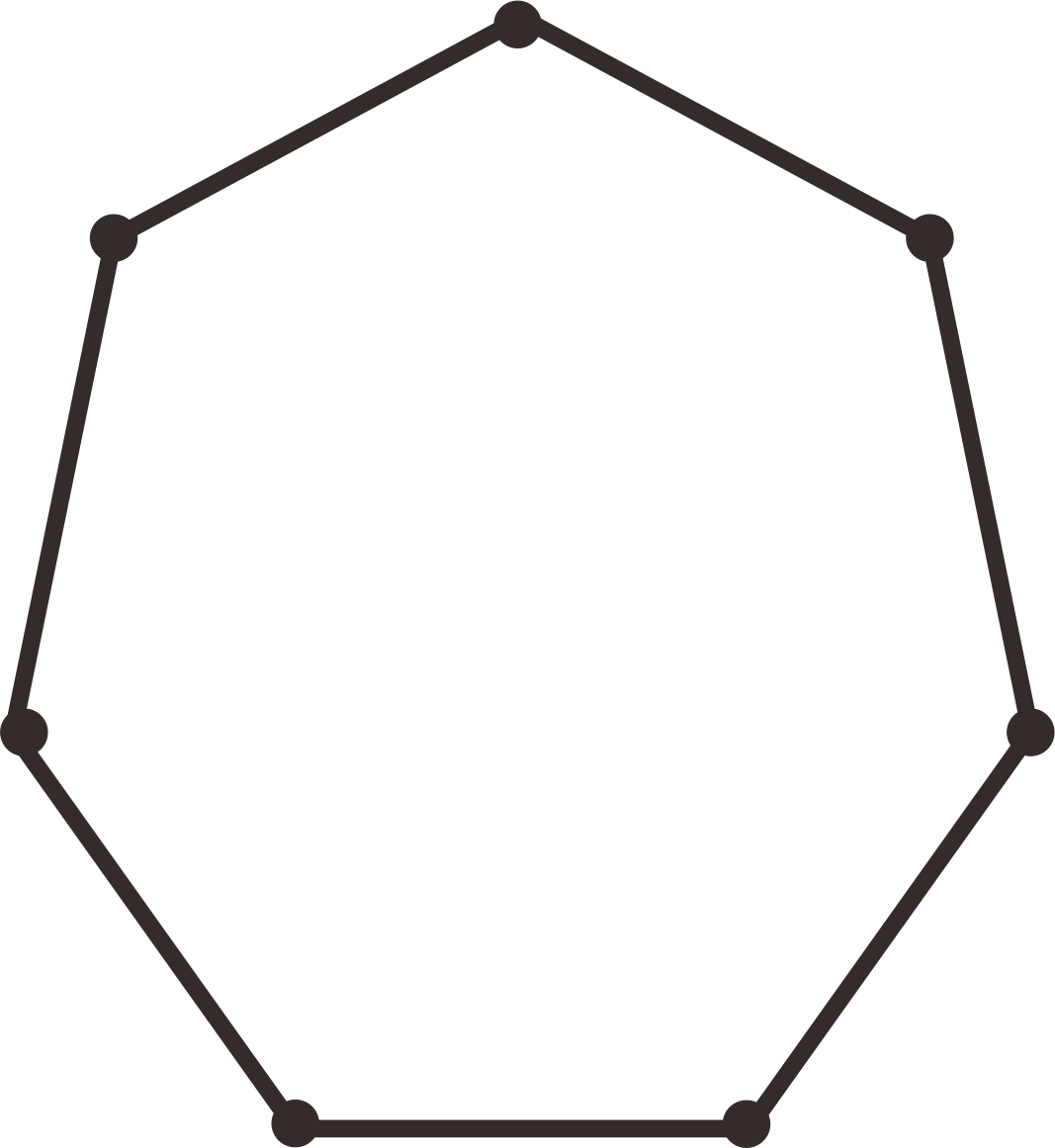}}
	\end{minipage}
    & \ \ \ \ 2
    & \begin{minipage}[b]{0.3\columnwidth}
		\centering
		\raisebox{-.25\height}{\includegraphics[width=0.7cm]{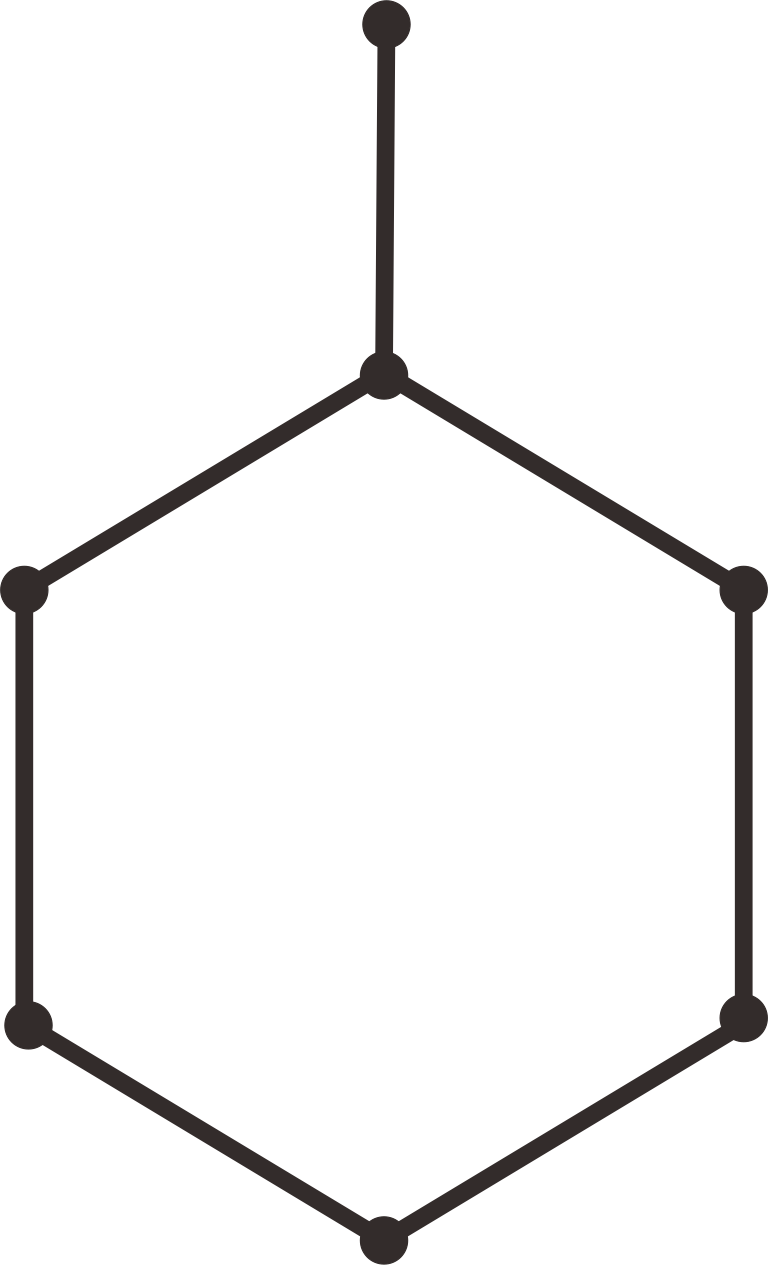}}
	\end{minipage}
    &  2.1602
    \\[2.5pt]
    \hline
     \begin{minipage}[b]{0.3\columnwidth}
		\centering
		\raisebox{-.25\height}{\includegraphics[width=1.5cm]{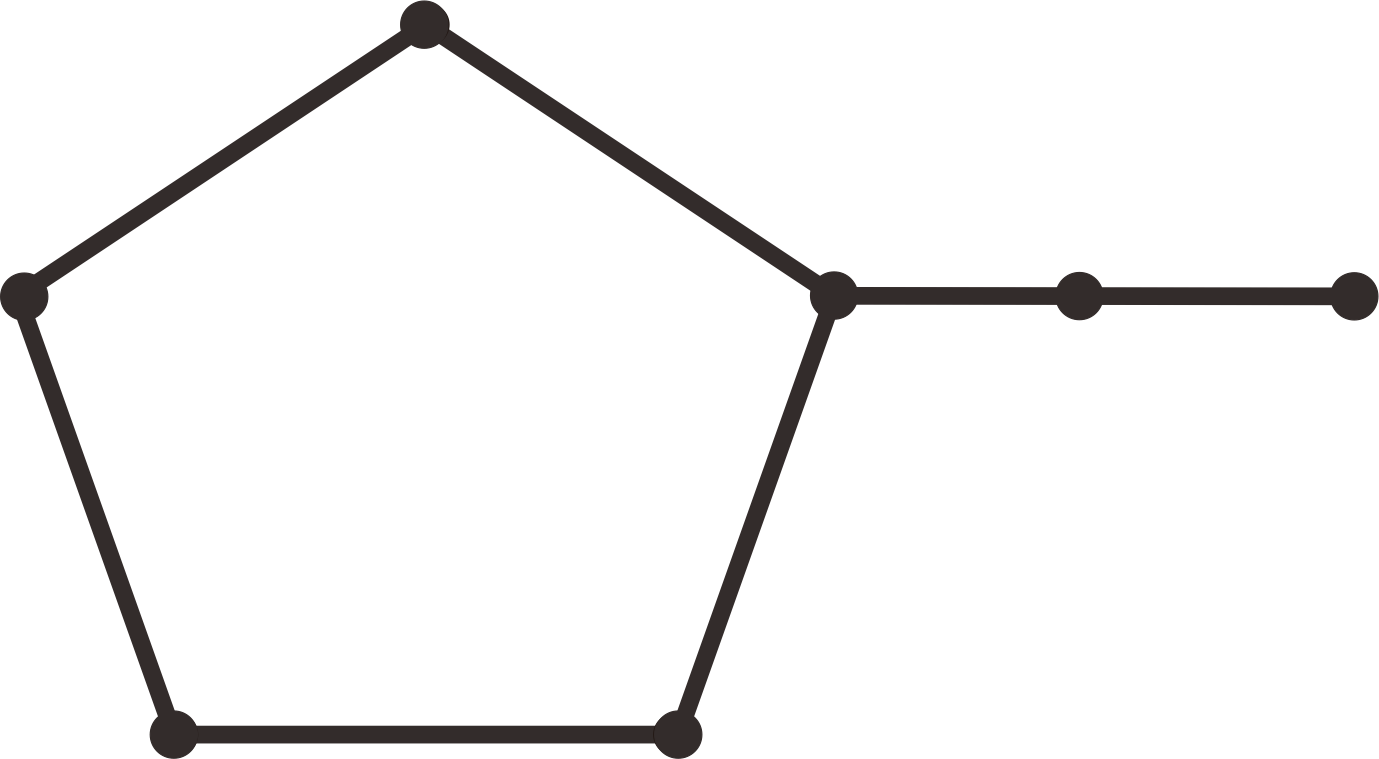}}
	\end{minipage}
    & 2.1827
    & \begin{minipage}[b]{0.3\columnwidth}
		\centering
		\raisebox{-.25\height}{\includegraphics[width=1.5cm]{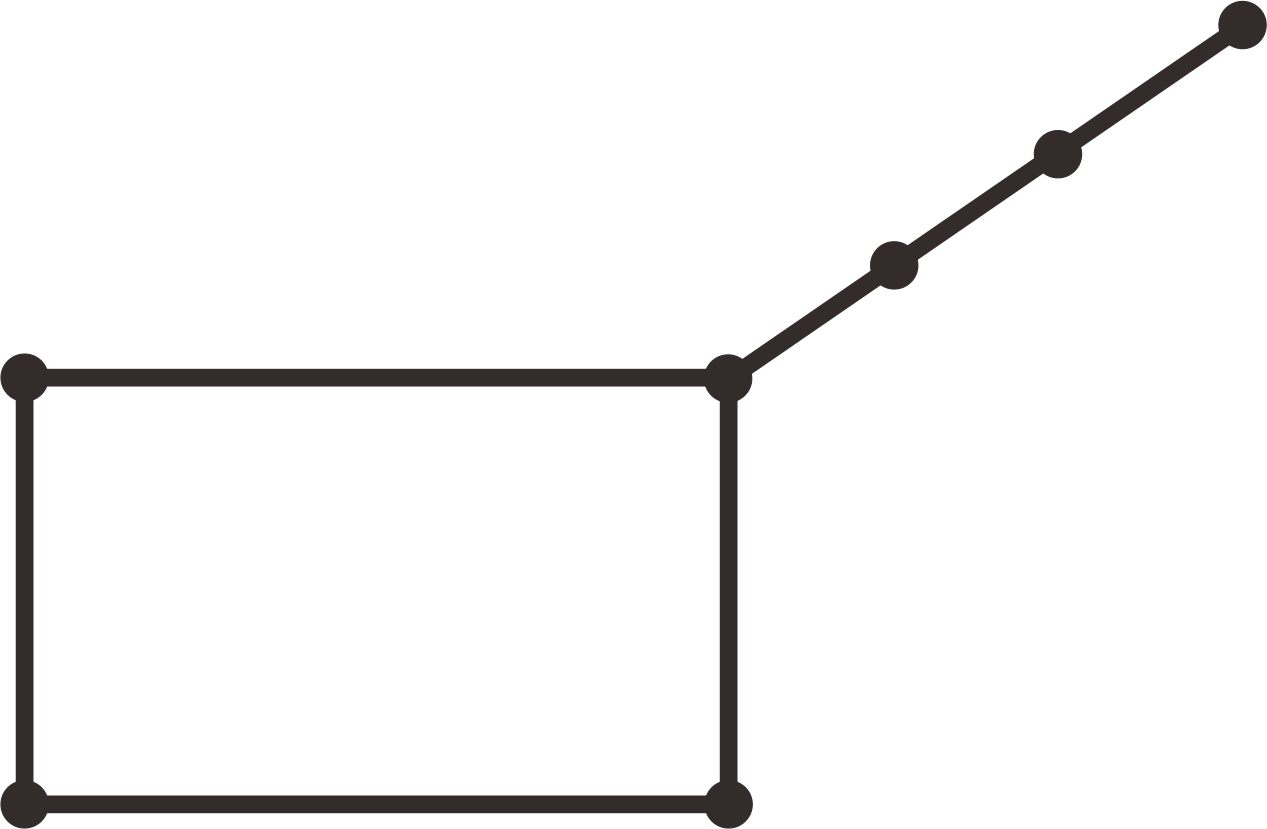}}
	\end{minipage}
    &  2.2188
    \\[2.5pt]
    \hline
    \begin{minipage}[b]{0.3\columnwidth}
		\centering
		\raisebox{-.25\height}{\includegraphics[width=1.6cm]{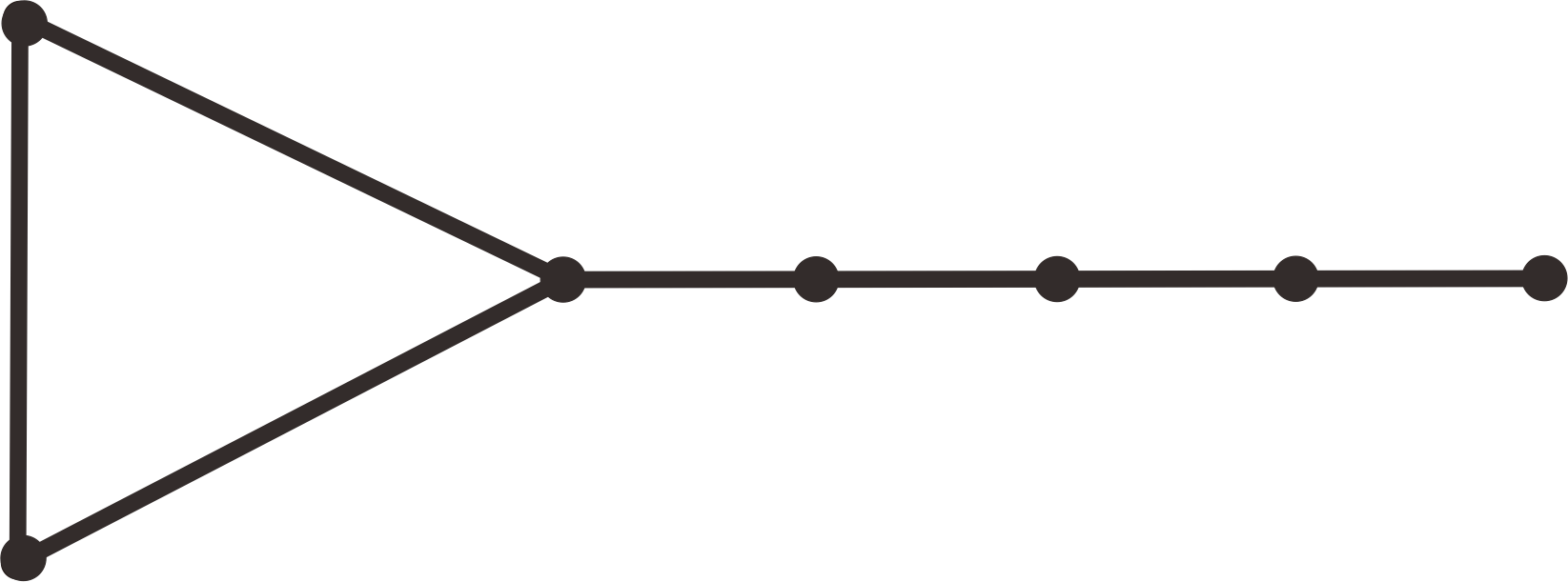}}
	\end{minipage}
    & 2.2661
    & \begin{minipage}[b]{0.3\columnwidth}
		\centering
		\raisebox{-.25\height}{\includegraphics[width=1.5cm]{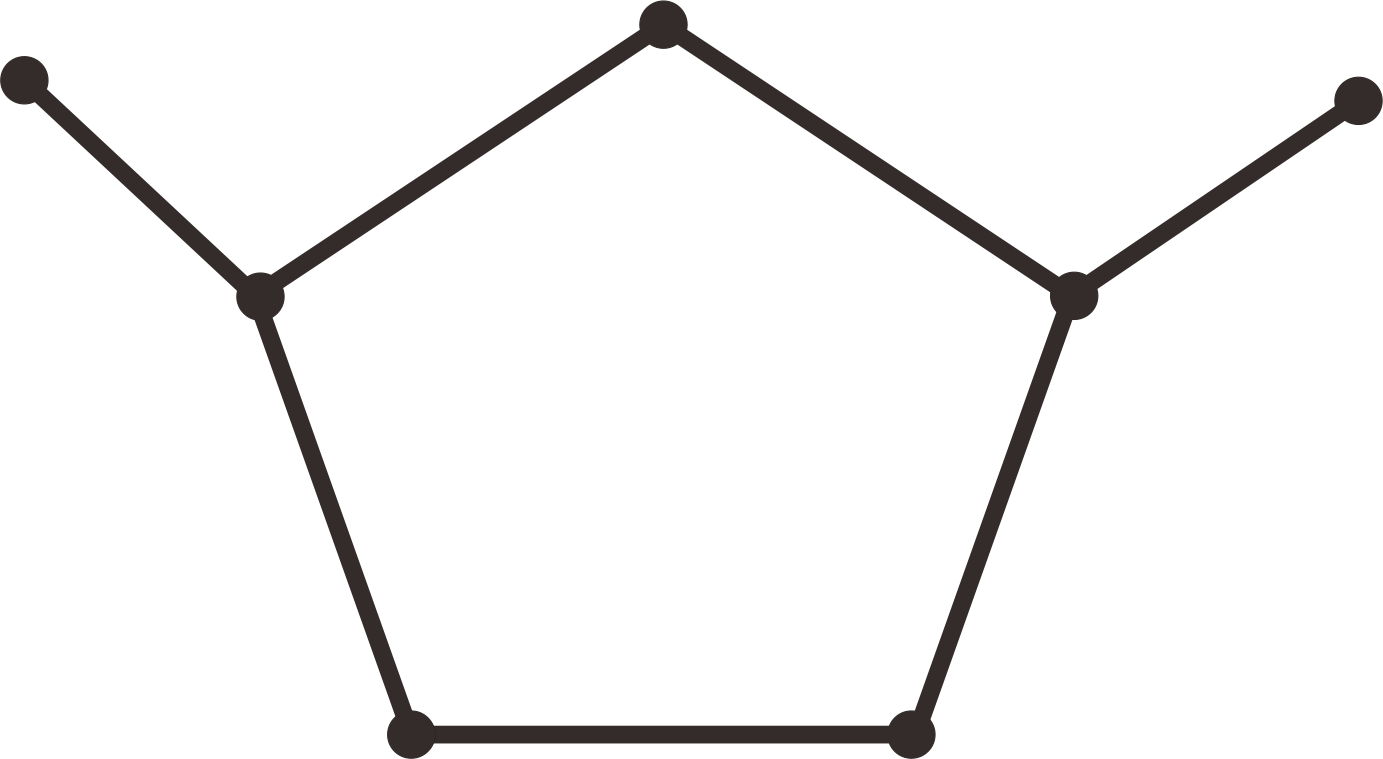}}
	\end{minipage}
    &  2.2942
    \\[2.5pt]
    \hline
     \begin{minipage}[b]{0.3\columnwidth}
		\centering
		\raisebox{-.30\height}{\includegraphics[width=1.3cm]{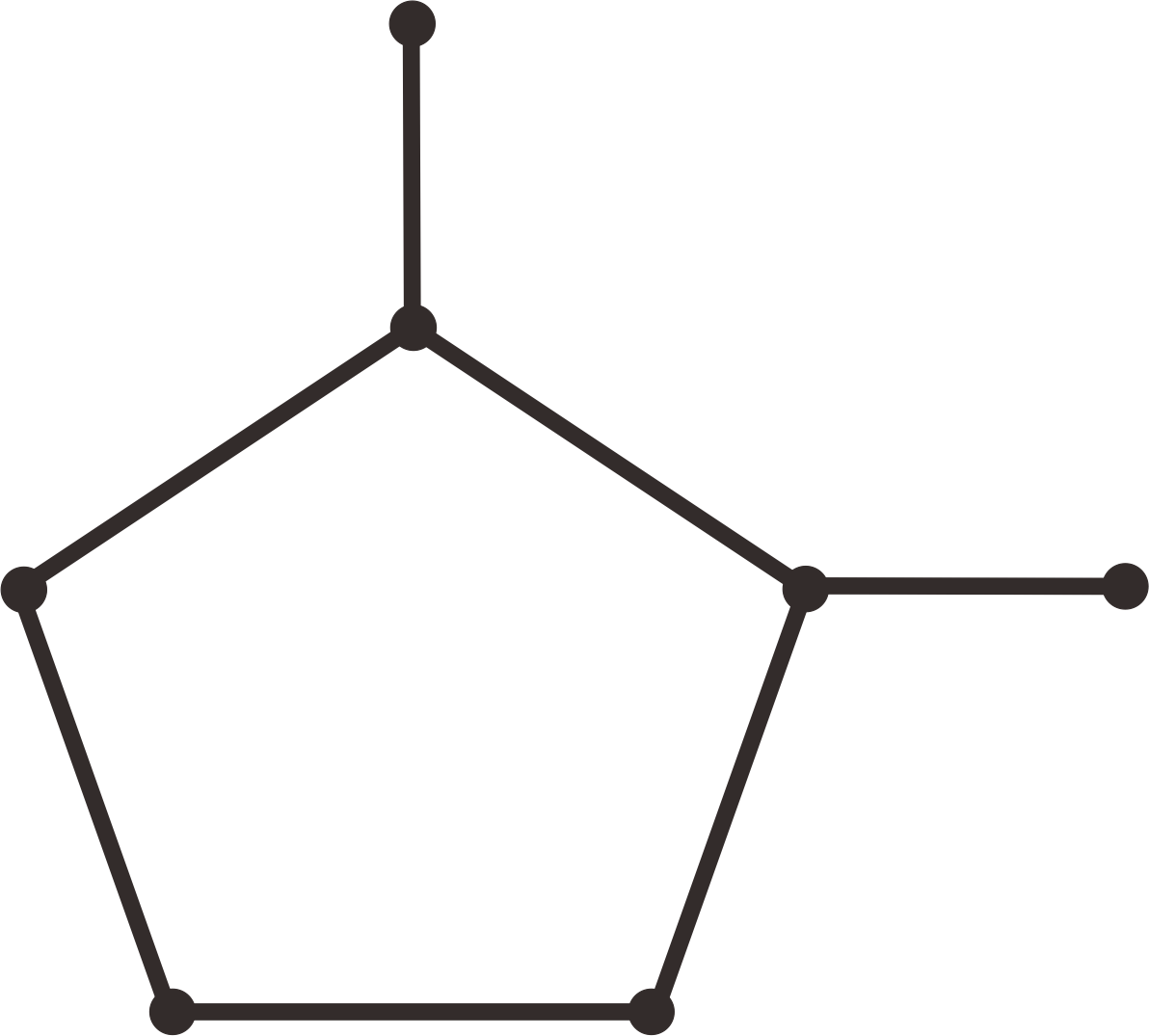}}
	\end{minipage}
    & 2.3041
    & \begin{minipage}[b]{0.3\columnwidth}
		\centering
		\raisebox{-.25\height}{\includegraphics[width=1.3cm]{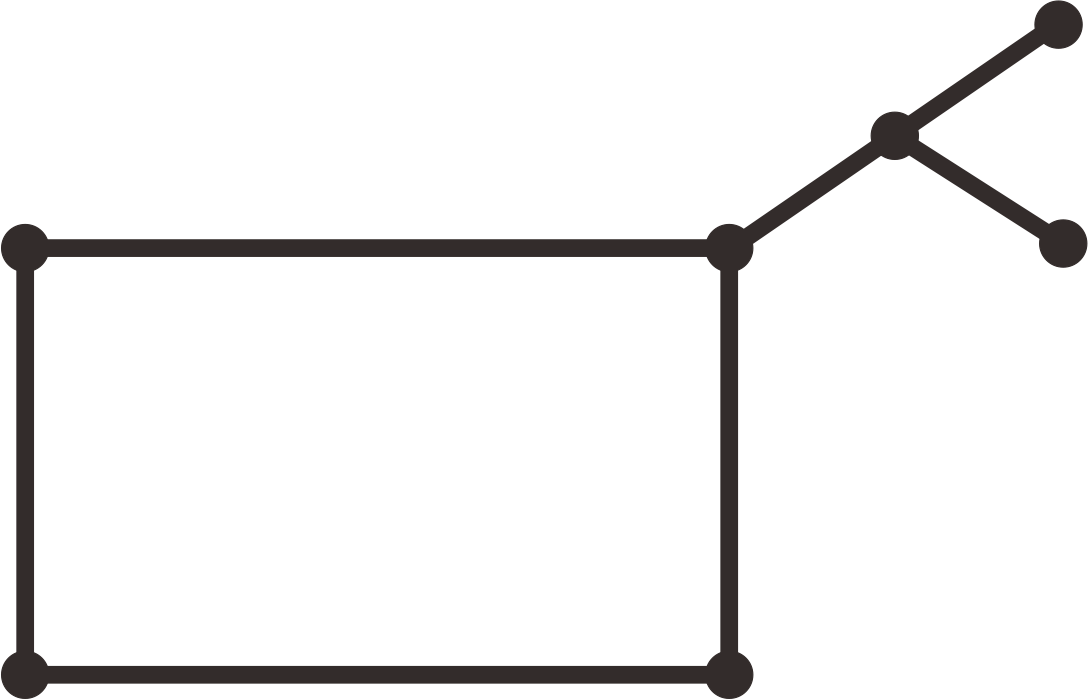}}
	\end{minipage}
    &  2.3094
    \\[2.5pt]
    \hline
     \begin{minipage}[b]{0.3\columnwidth}
		\centering
		\raisebox{-.25\height}{\includegraphics[width=1.7cm]{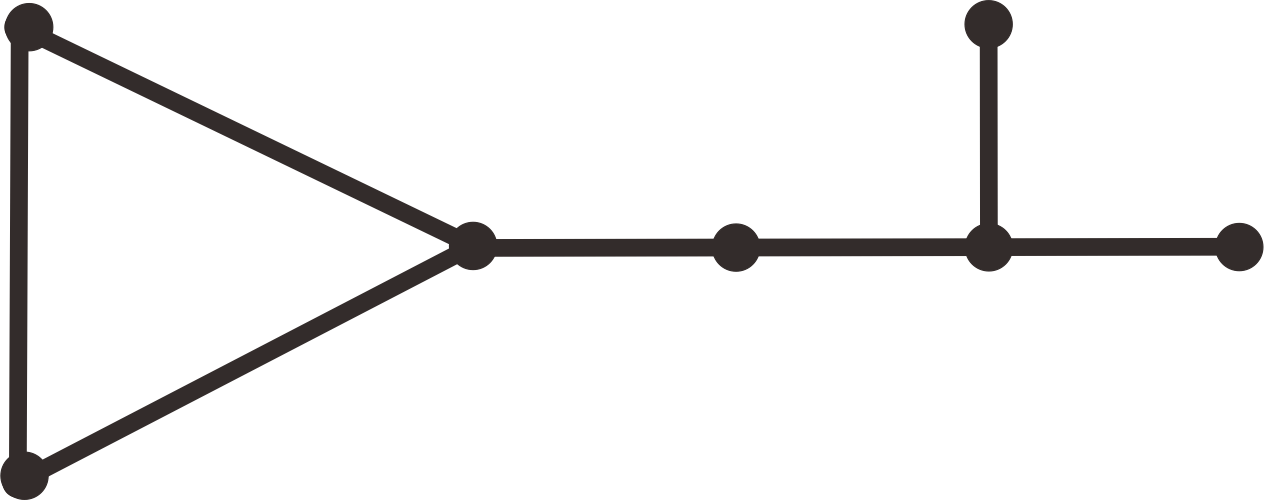}}
	\end{minipage}
    & 2.316
    & \begin{minipage}[b]{0.3\columnwidth}
		\centering
		\raisebox{-.25\height}{\includegraphics[width=1.3cm]{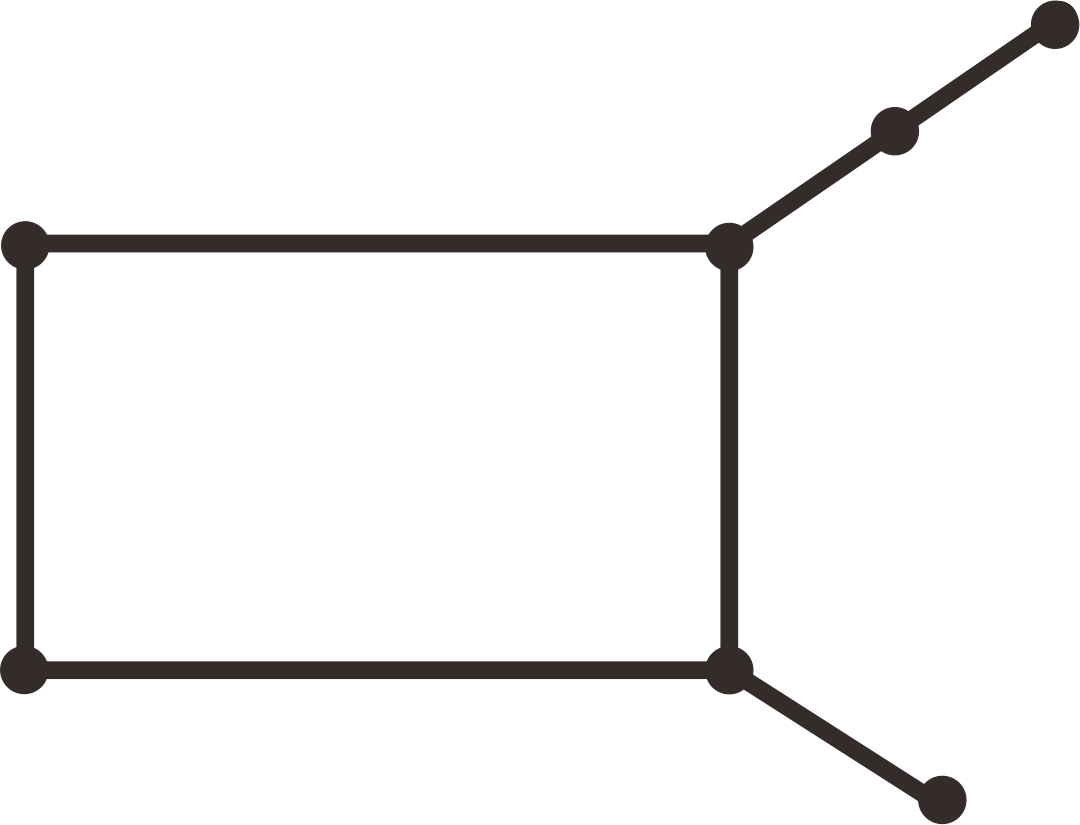}}
	\end{minipage}
    &  2.3413
    \\[2.5pt]
    \hline
    \begin{minipage}[b]{0.3\columnwidth}
		\centering
		\raisebox{-.25\height}{\includegraphics[width=1.5cm]{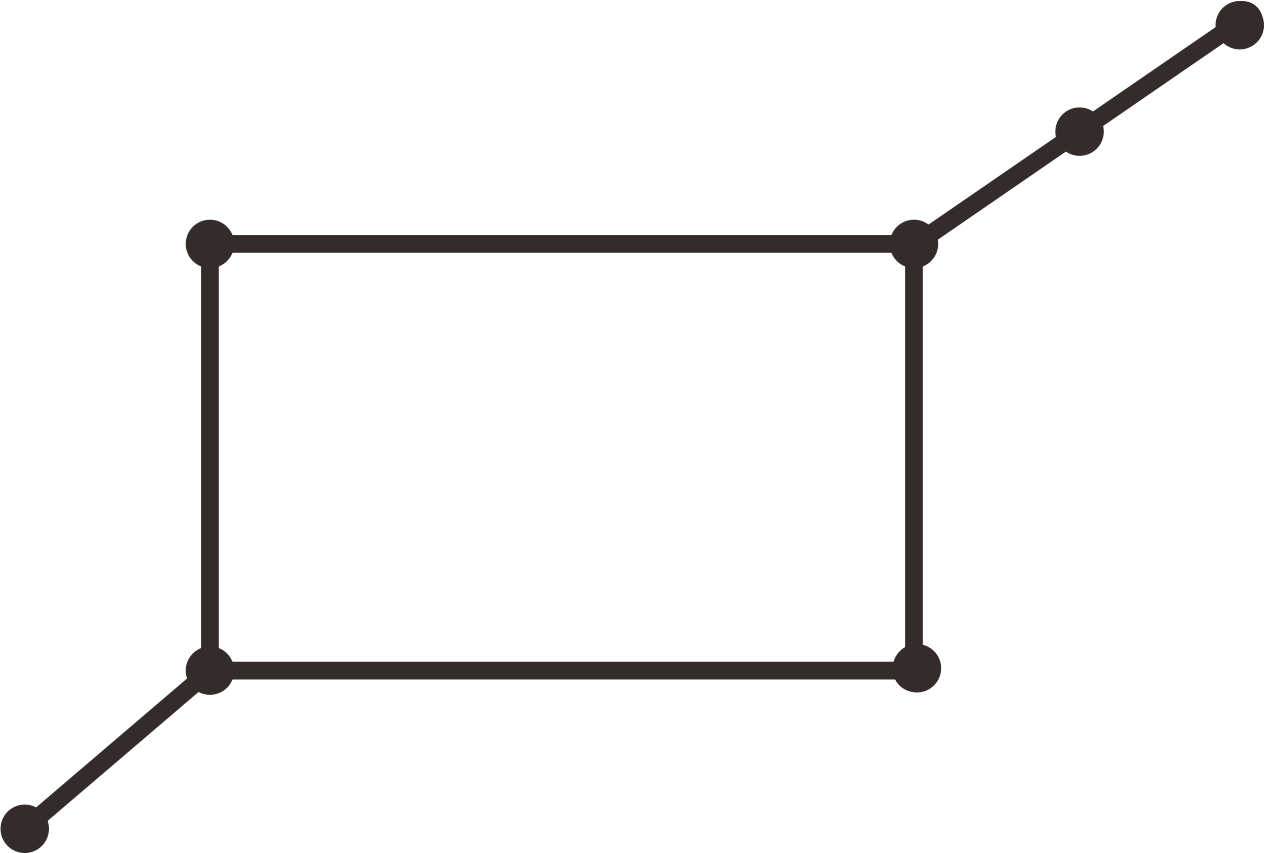}}
	\end{minipage}
    & 2.3428
    & \begin{minipage}[b]{0.3\columnwidth}
		\centering
		\raisebox{-.25\height}{\includegraphics[width=1.7cm]{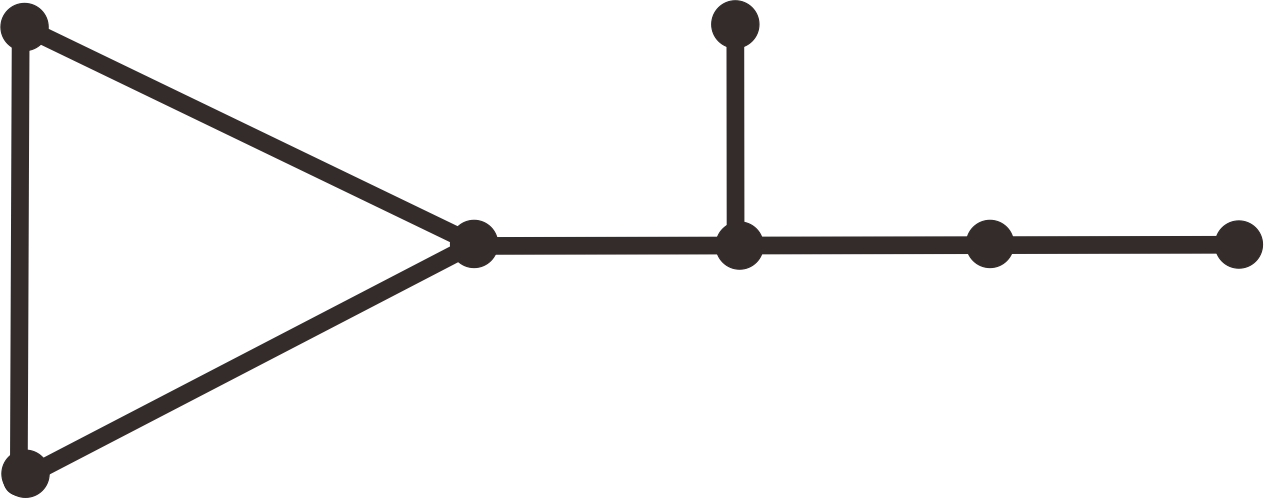}}
	\end{minipage}
    &  2.3528
    \\[2.5pt]
    \hline
    \begin{minipage}[b]{0.3\columnwidth}
		\centering
		\raisebox{-.30\height}{\includegraphics[width=1.3cm]{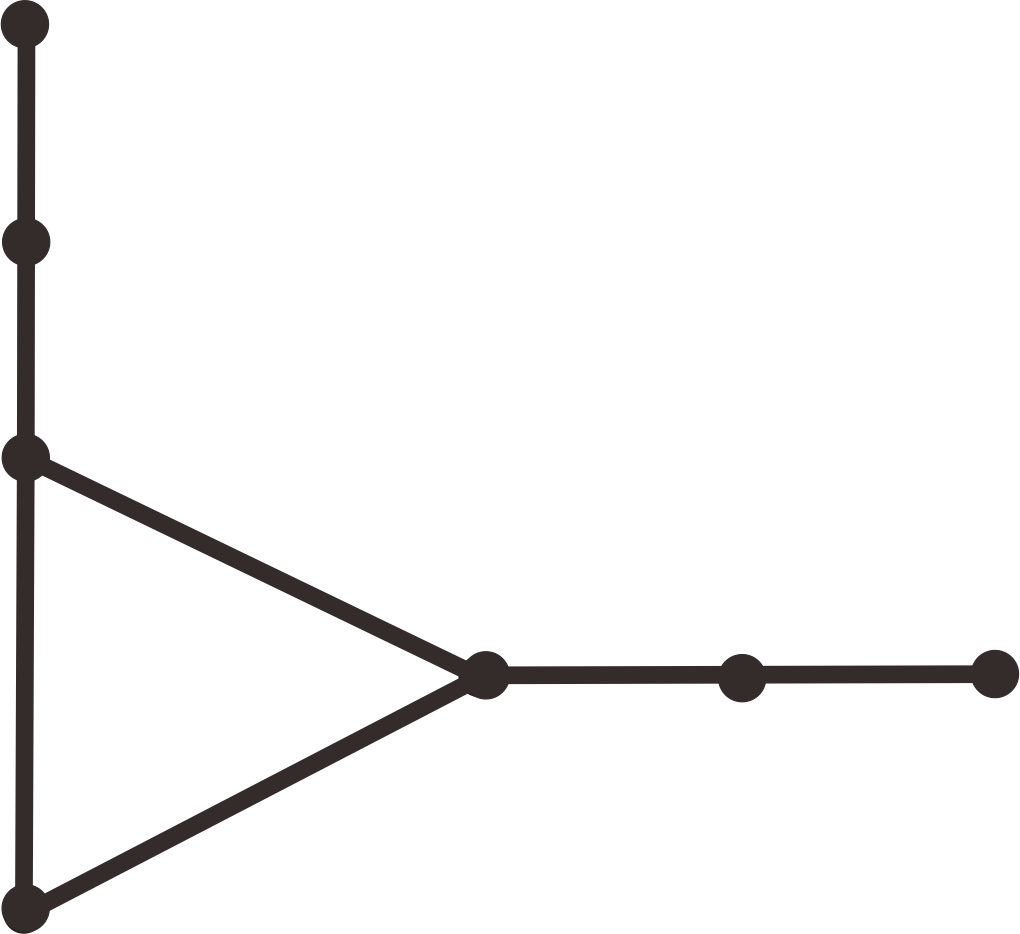}}
	\end{minipage}
    & 2.4044
    & \begin{minipage}[b]{0.3\columnwidth}
		\centering
		\raisebox{-.25\height}{\includegraphics[width=1.6cm]{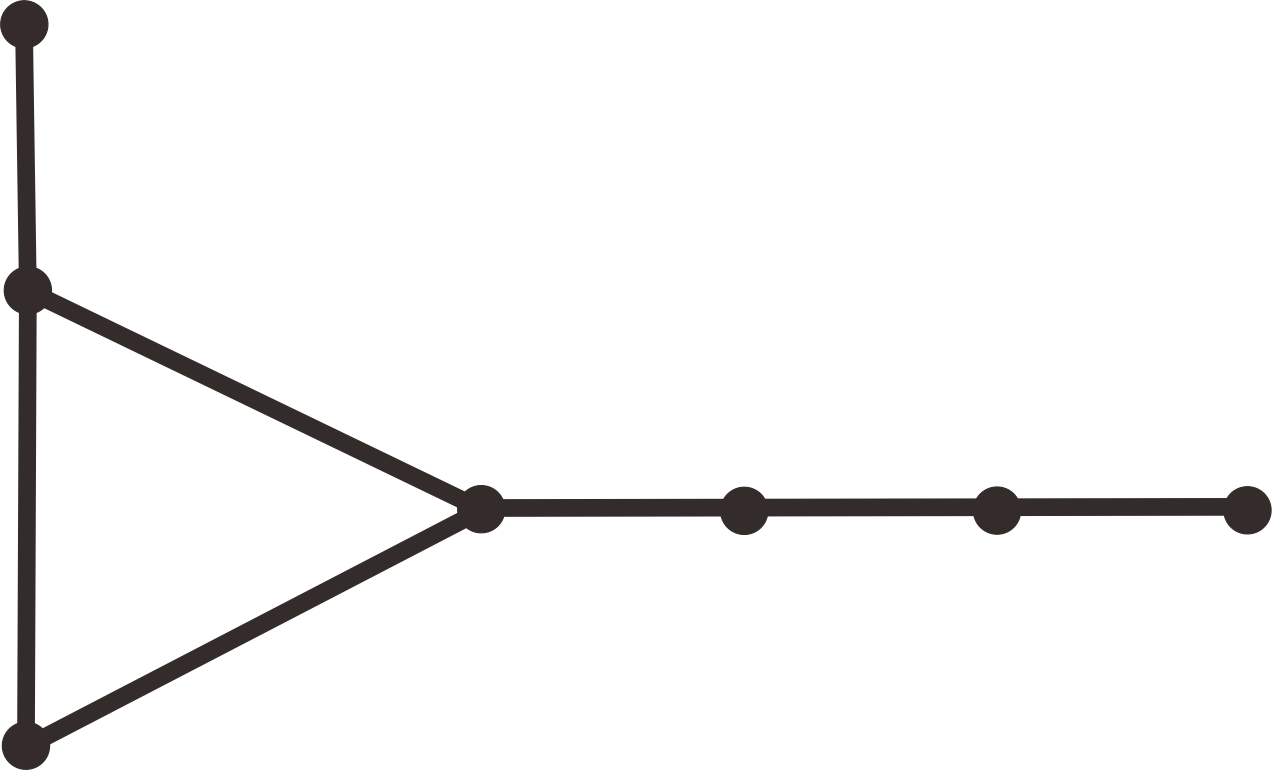}}
	\end{minipage}
    &  2.4135
    \\[2.5pt]
   \hline
   \begin{minipage}[b]{0.3\columnwidth}
		\centering
		\raisebox{-.25\height}{\includegraphics[width=1.3cm]{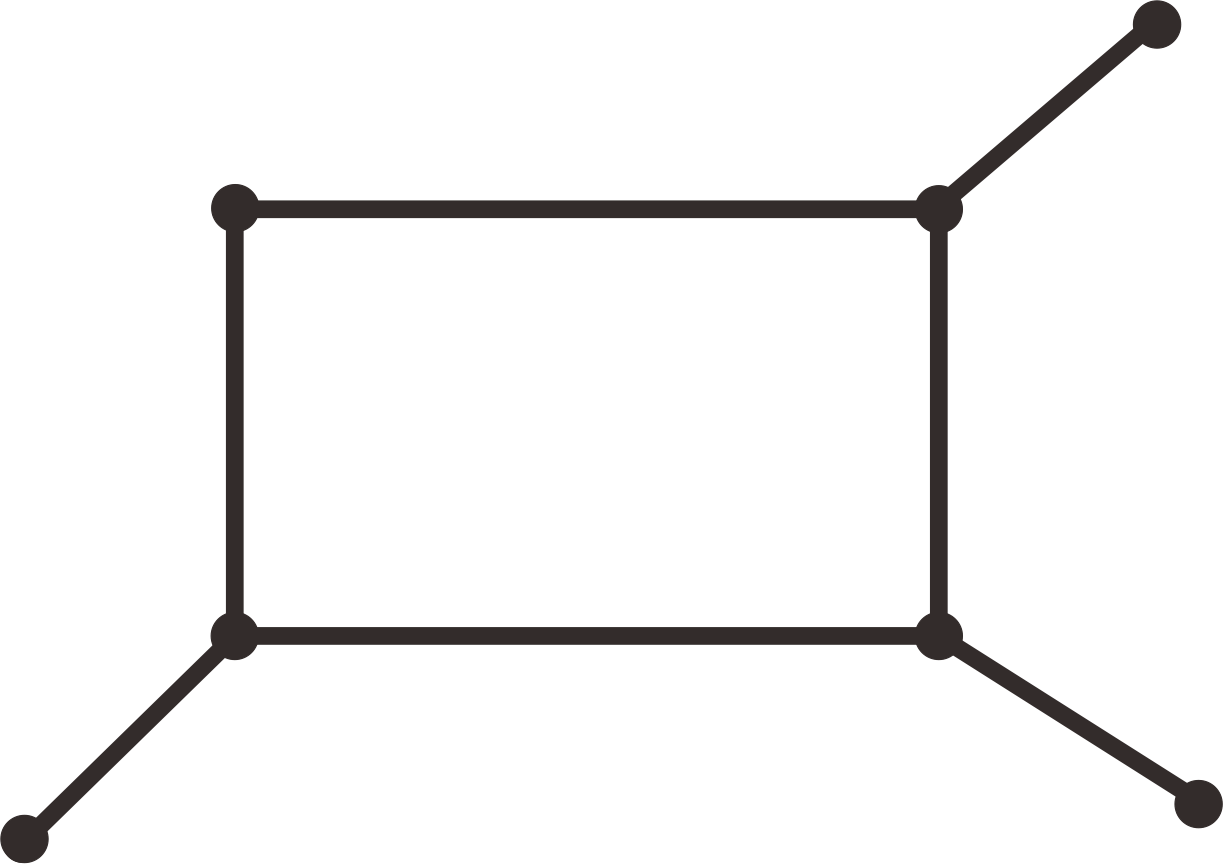}}
	\end{minipage}
    & 2.4473
    & \begin{minipage}[b]{0.3\columnwidth}
		\centering
		\raisebox{-.25\height}{\includegraphics[width=1.3cm]{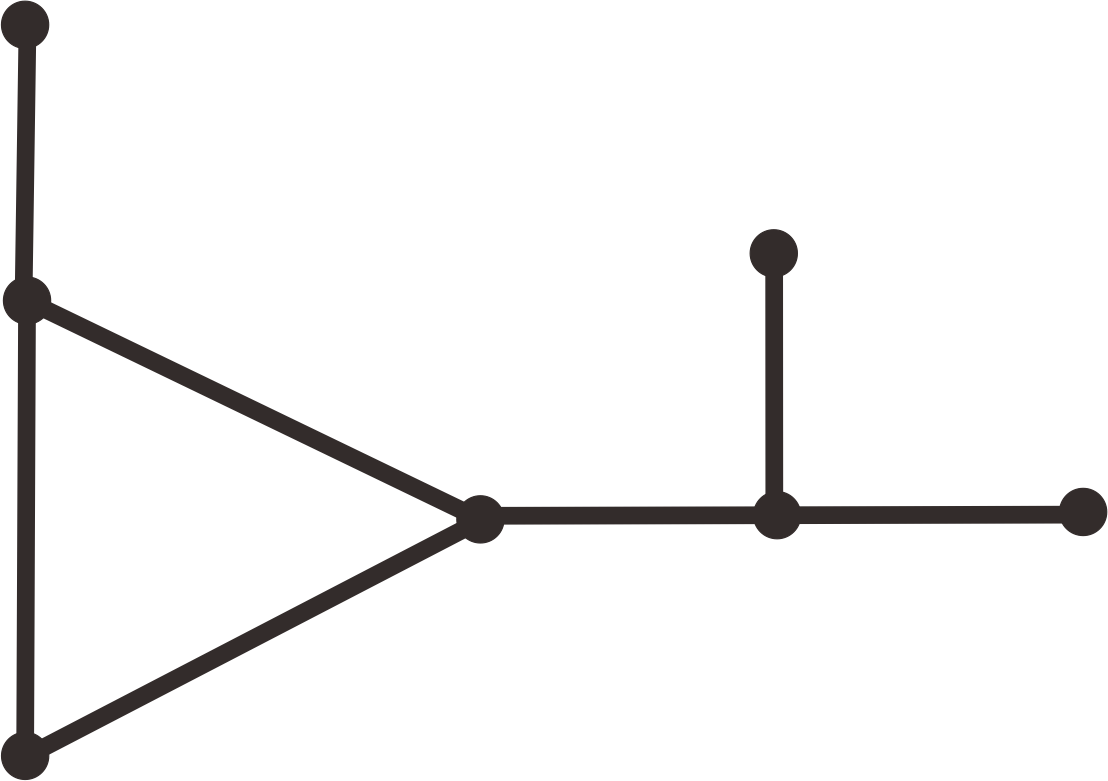}}
	\end{minipage}
    &  2.4732
    \\[2.5pt]
    \hline
    \begin{minipage}[b]{0.3\columnwidth}
		\centering
		\raisebox{-.25\height}{\includegraphics[width=1.3cm]{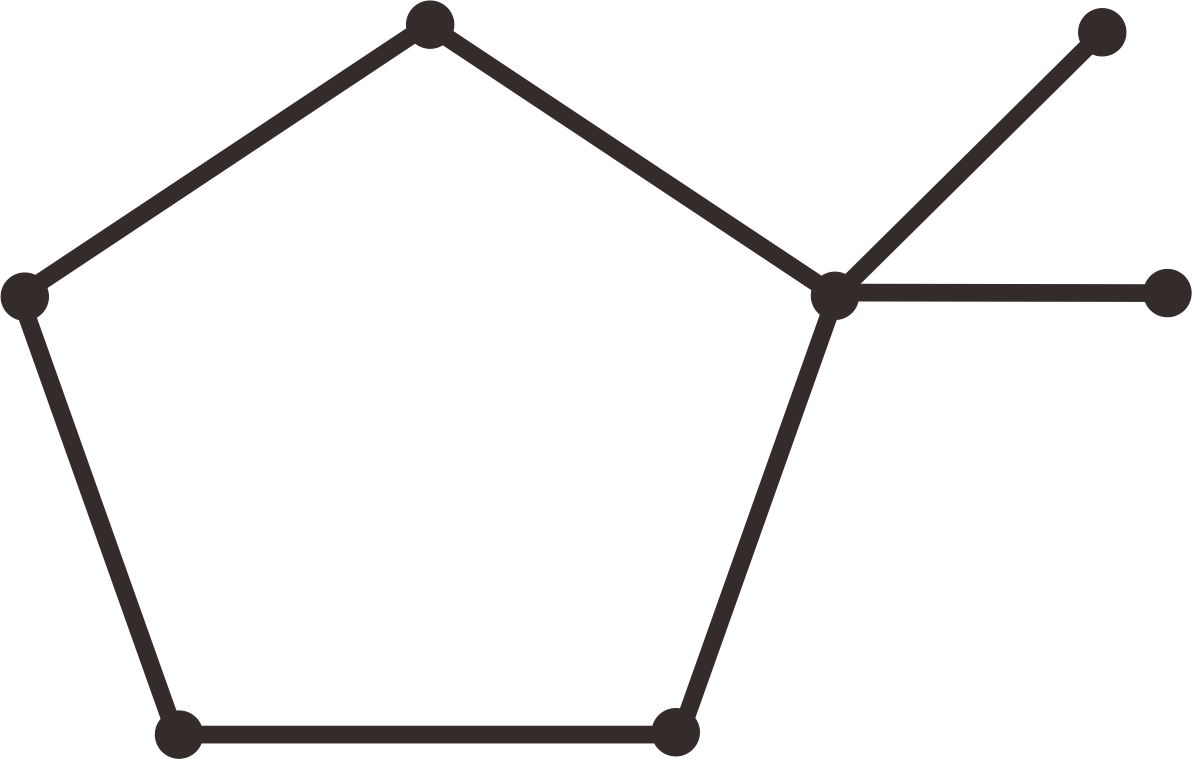}}
	\end{minipage}
    & 2.4908
    & \begin{minipage}[b]{0.3\columnwidth}
		\centering
		\raisebox{-.25\height}{\includegraphics[width=1.3cm]{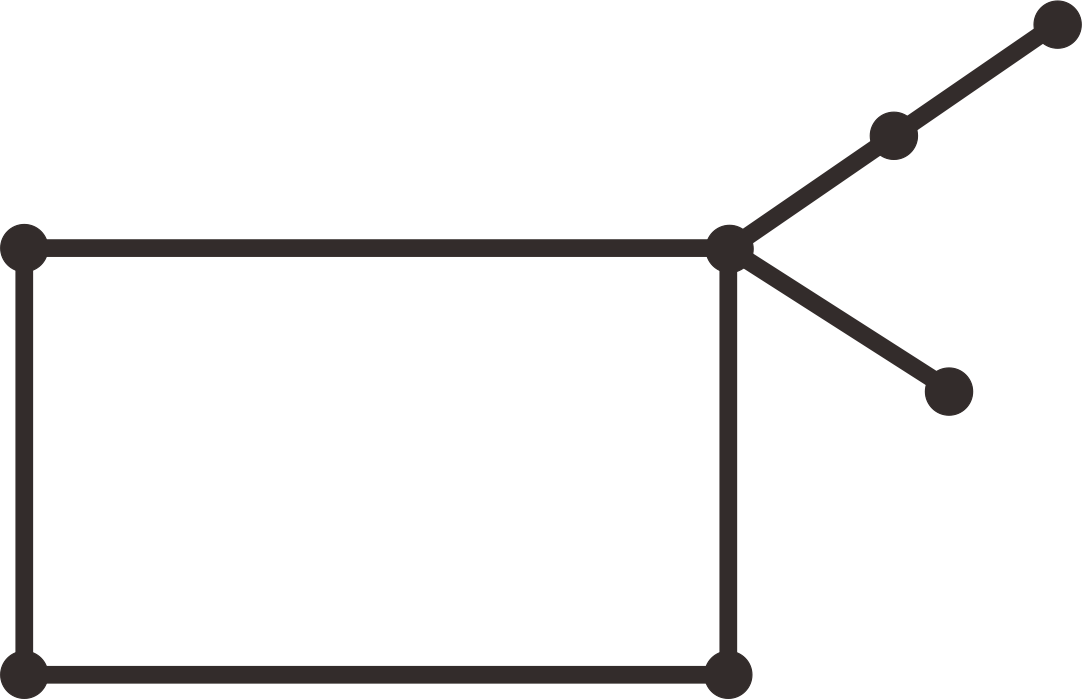}}
	\end{minipage}
    &  2.499
    \\[2.5pt]
   \hline
   \begin{minipage}[b]{0.3\columnwidth}
		\centering
		\raisebox{-.30\height}{\includegraphics[width=1cm]{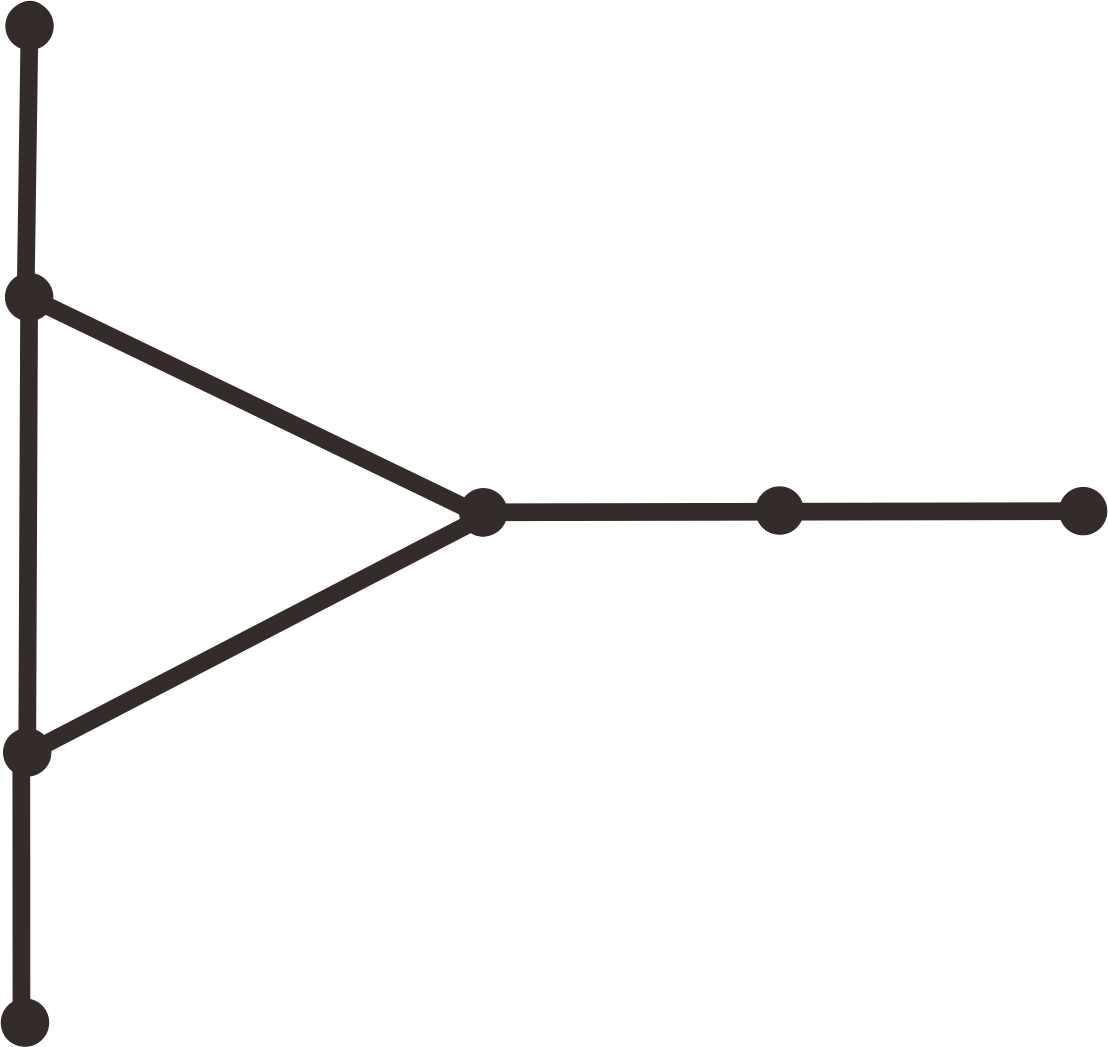}}
	\end{minipage}
    & 2.5202
    & \begin{minipage}[b]{0.3\columnwidth}
		\centering
		\raisebox{-.25\height}{\includegraphics[width=1.3cm]{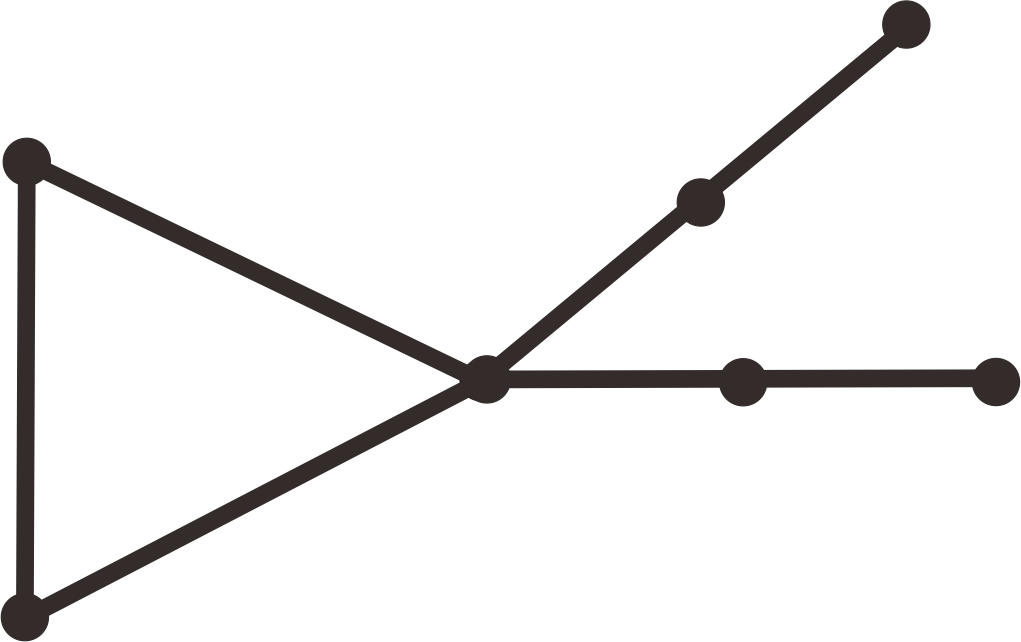}}
	\end{minipage}
    &  2.5376
    \\[2.5pt]
    \hline
     \begin{minipage}[b]{0.3\columnwidth}
		\centering
		\raisebox{-.25\height}{\includegraphics[width=1.7cm]{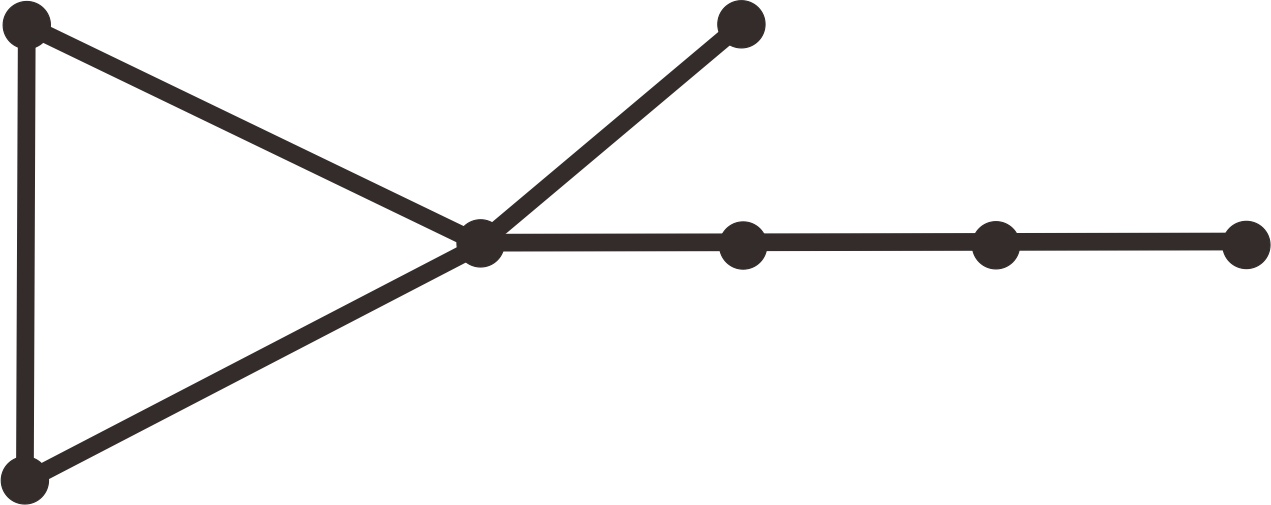}}
	\end{minipage}
    & 2.5727
    & \begin{minipage}[b]{0.3\columnwidth}
		\centering
		\raisebox{-.25\height}{\includegraphics[width=1.3cm]{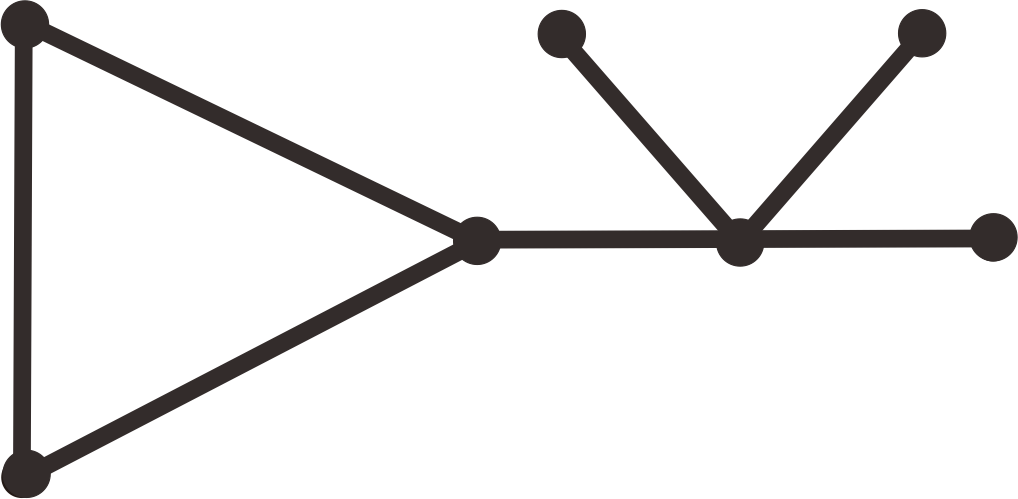}}
	\end{minipage}
    &  2.5962
    \\[2.5pt]
    \hline
    \begin{minipage}[b]{0.3\columnwidth}
		\centering
		\raisebox{-.25\height}{\includegraphics[width=1.3cm]{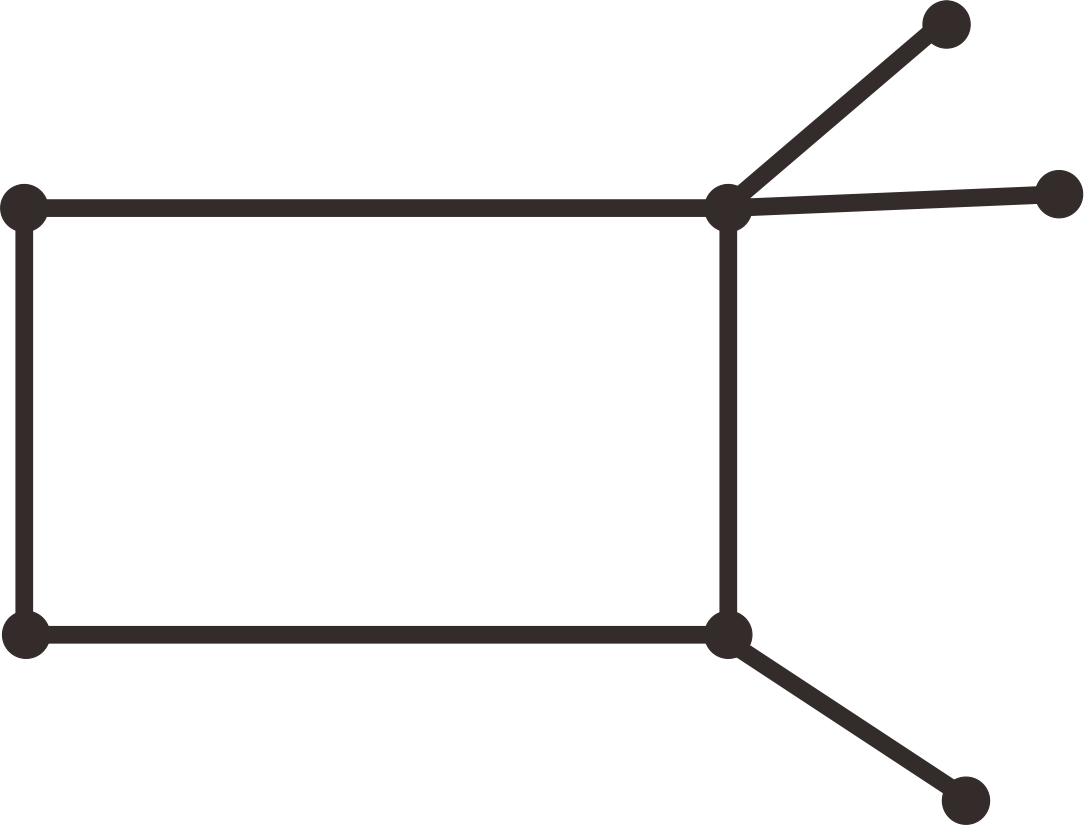}}
	\end{minipage}
    & 2.5992
    & \begin{minipage}[b]{0.3\columnwidth}
		\centering
		\raisebox{-.25\height}{\includegraphics[width=1.6cm]{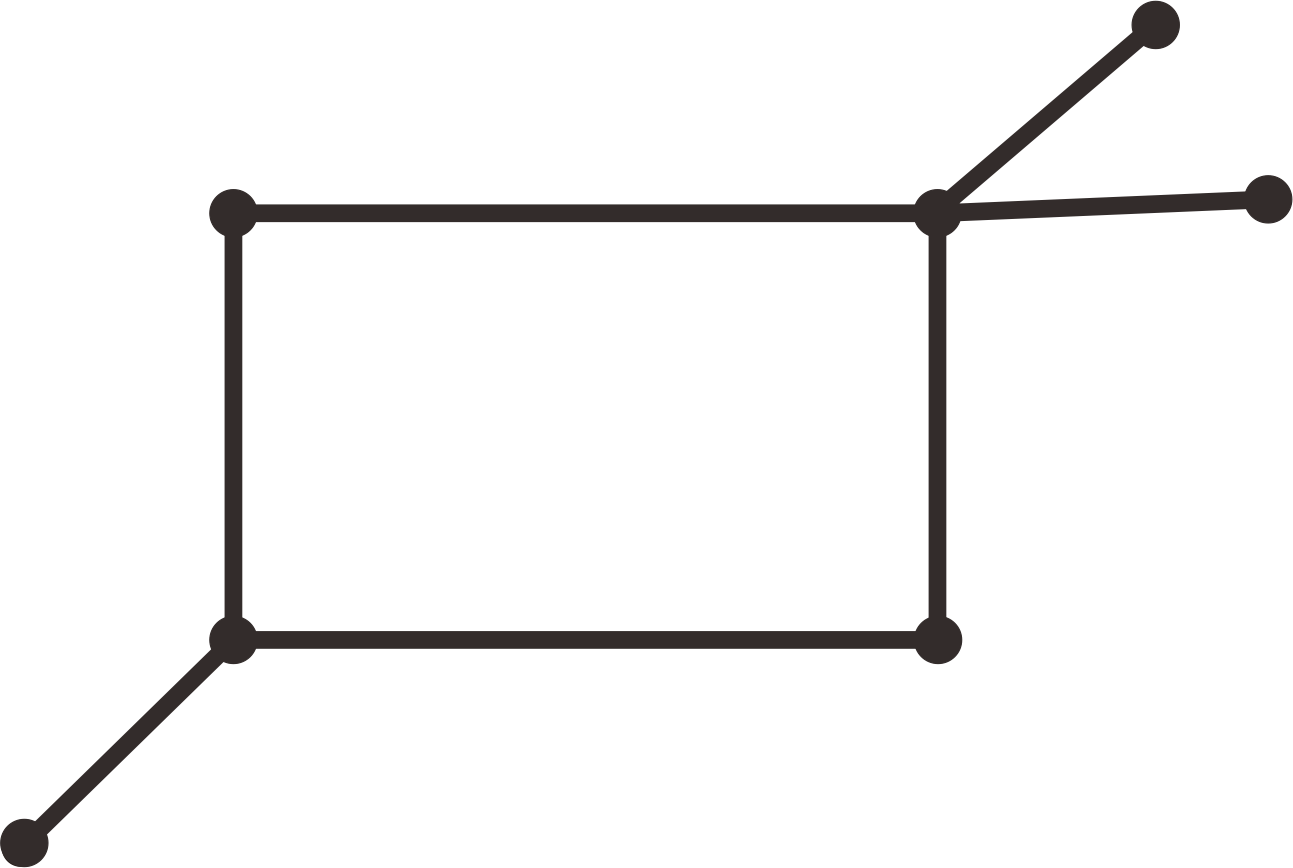}}
	\end{minipage}
    &  2.6023
    \\[2.5pt]
    \hline
     \begin{minipage}[b]{0.3\columnwidth}
		\centering
		\raisebox{-.25\height}{\includegraphics[width=1.3cm]{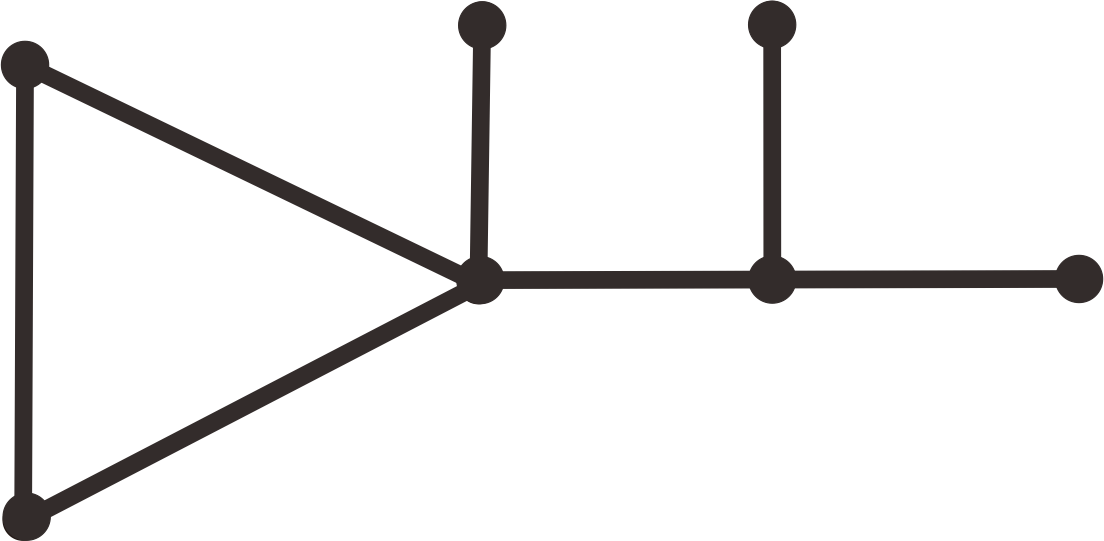}}
	\end{minipage}
    & 2.6209
    & \begin{minipage}[b]{0.3\columnwidth}
		\centering
		\raisebox{-.25\height}{\includegraphics[width=1.3cm]{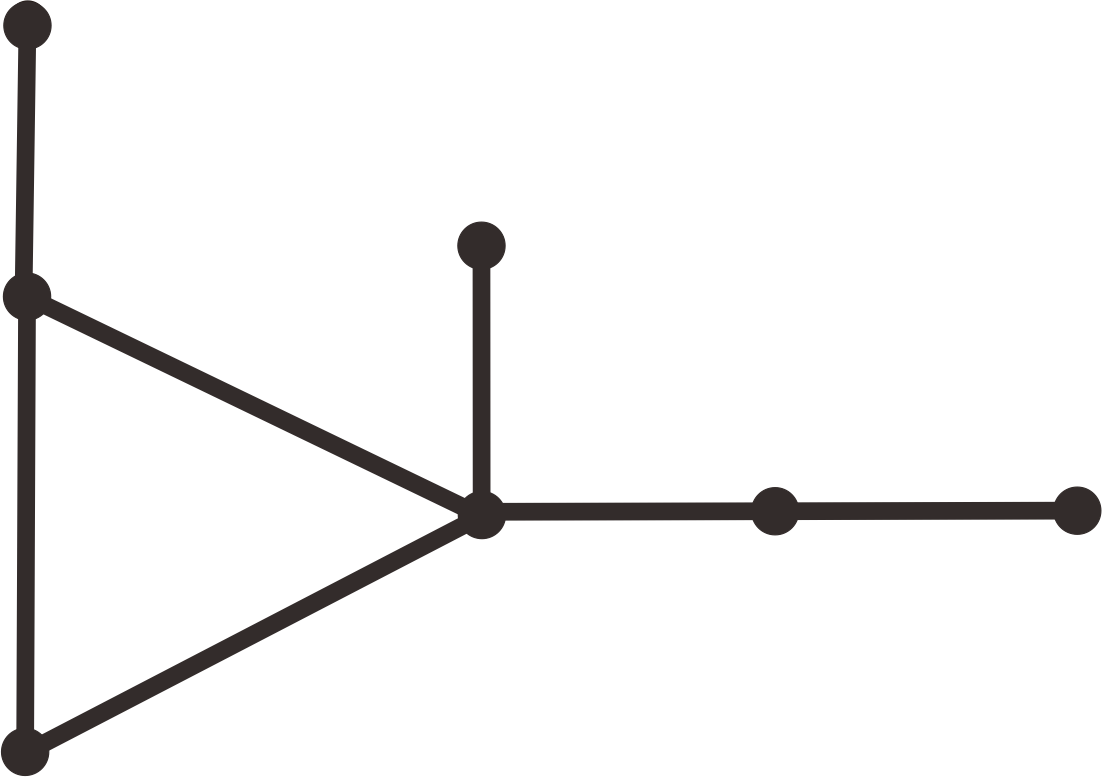}}
	\end{minipage}
    &  2.6564
    \\[2.5pt]
    \hline
     \begin{minipage}[b]{0.3\columnwidth}
		\centering
		\raisebox{-.25\height}{\includegraphics[width=1.5cm]{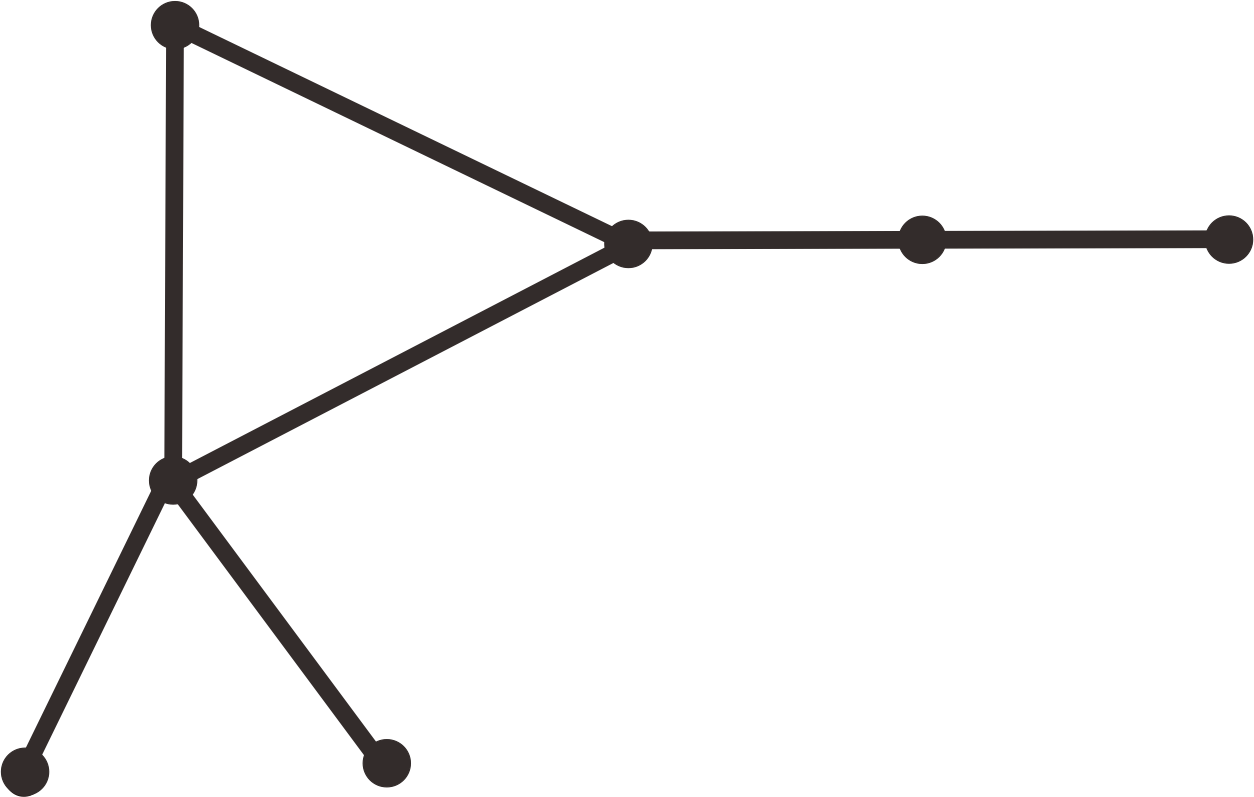}}
	\end{minipage}
    & 2.6795
    & \begin{minipage}[b]{0.3\columnwidth}
		\centering
		\raisebox{-.25\height}{\includegraphics[width=1cm]{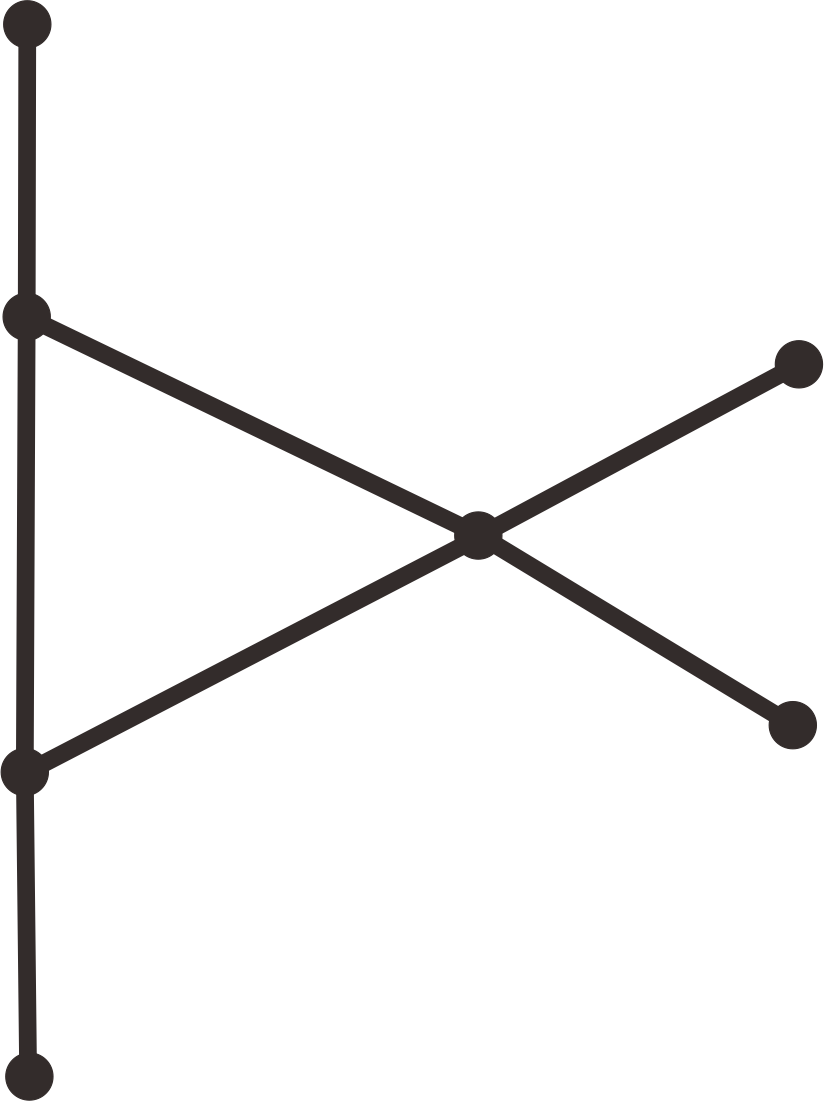}}
	\end{minipage}
    &  2.75
    \\[2.5pt]
    \hline
     \begin{minipage}[b]{0.3\columnwidth}
		\centering
		\raisebox{-.25\height}{\includegraphics[width=1.4cm]{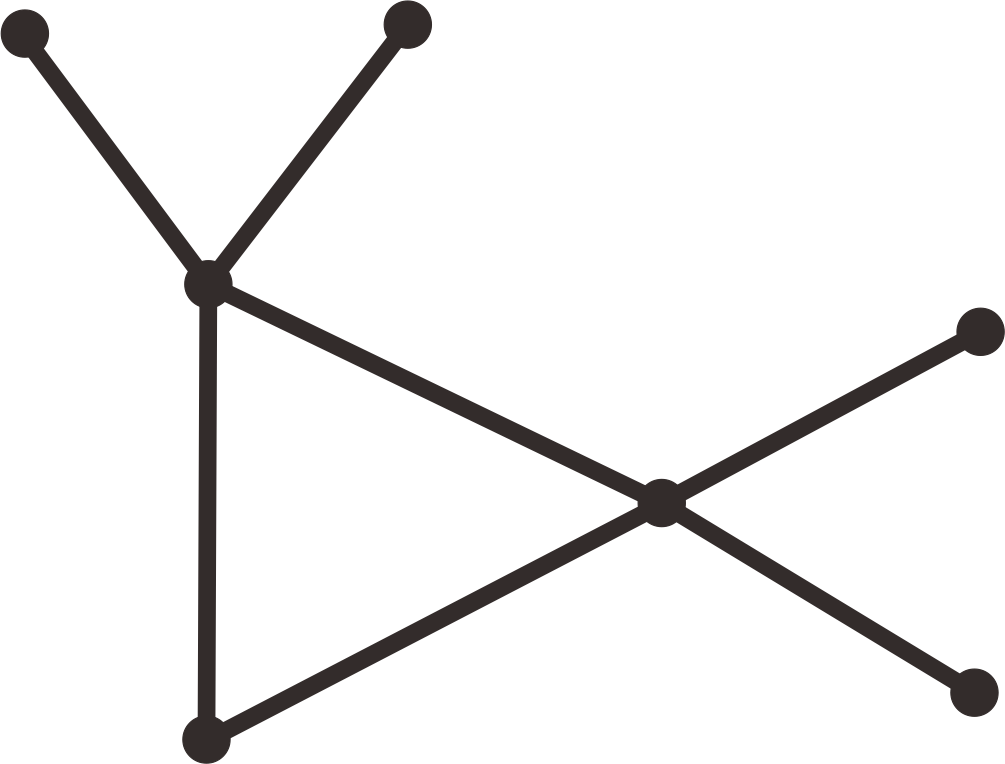}}
	\end{minipage}
    & 2.8717
    & \begin{minipage}[b]{0.3\columnwidth}
		\centering
		\raisebox{-.25\height}{\includegraphics[width=1.4cm]{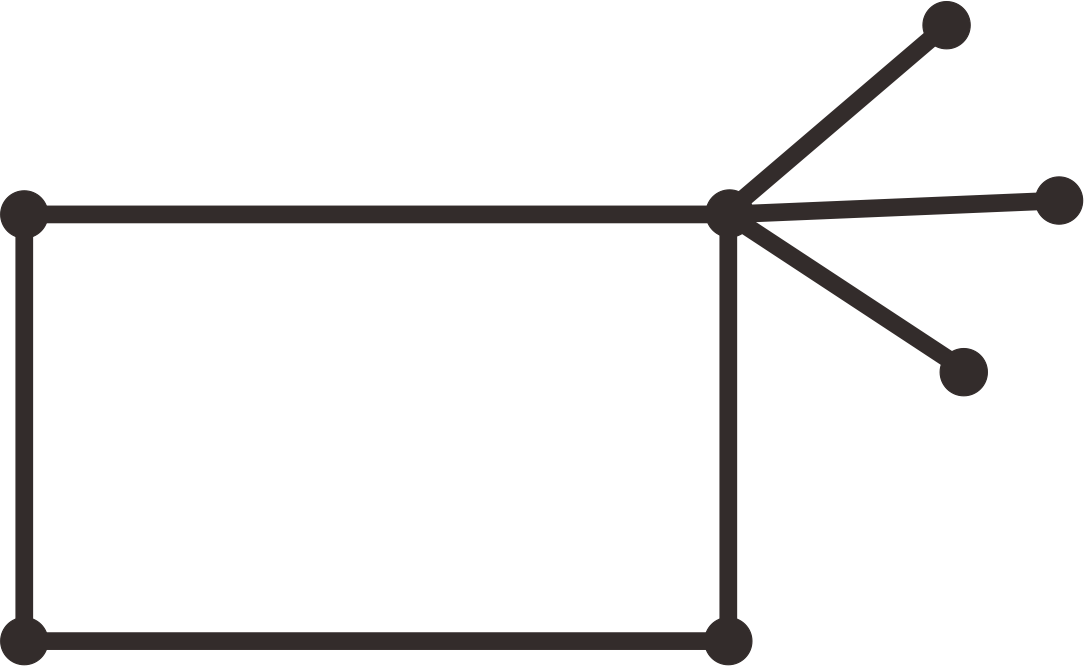}}
	\end{minipage}
    &  2.9314
    \\[2.5pt]
    \hline
     \begin{minipage}[b]{0.3\columnwidth}
		\centering
		\raisebox{-.25\height}{\includegraphics[width=1.3cm]{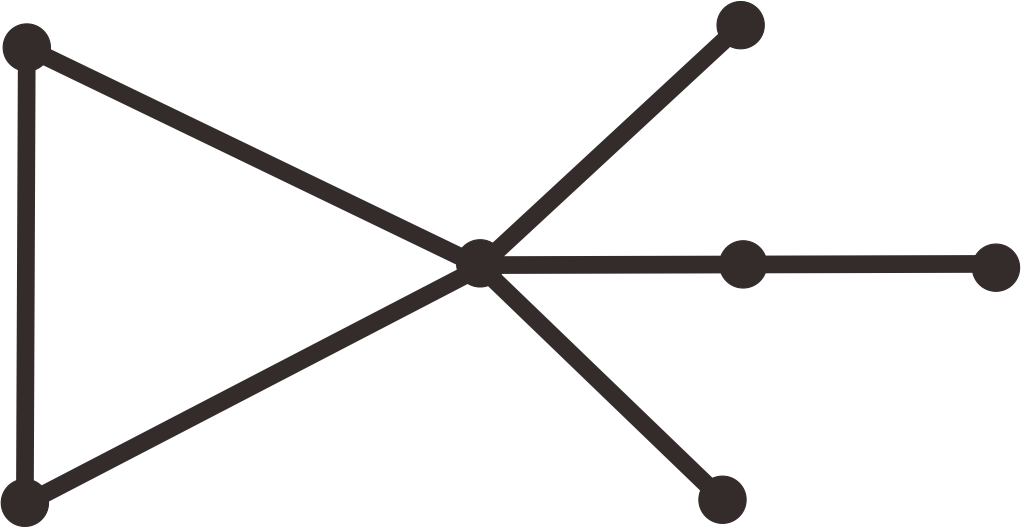}}
	\end{minipage}
    & 2.9516
    & \begin{minipage}[b]{0.3\columnwidth}
		\centering
		\raisebox{-.35\height}{\includegraphics[width=1.3cm]{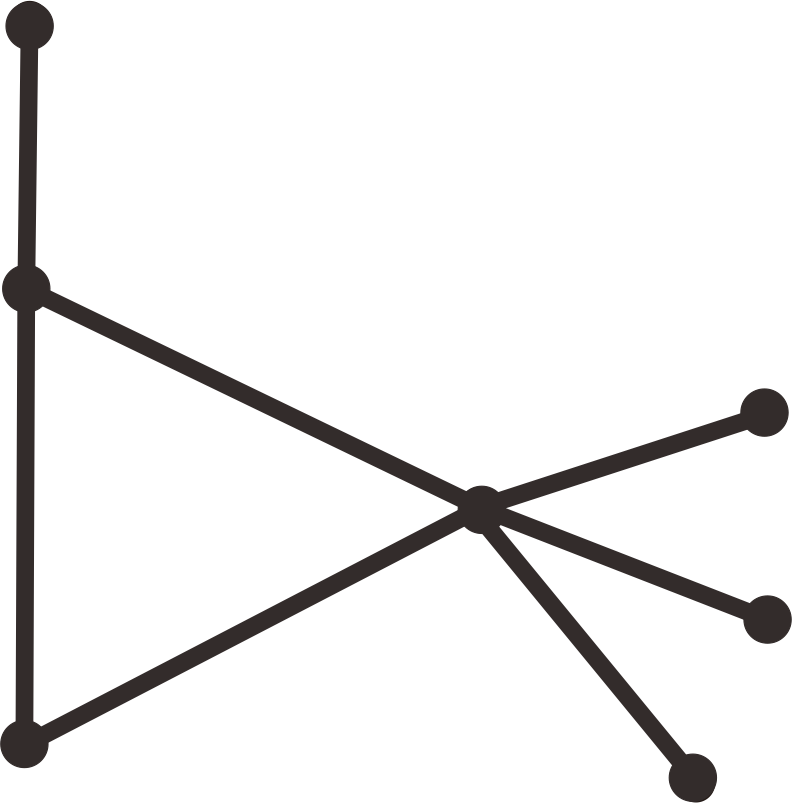}}
	\end{minipage}
    &  3.0453
    \\[2.5pt]
    \hline
    \begin{minipage}[b]{0.3\columnwidth}
		\centering
		\raisebox{-.25\height}{\includegraphics[width=1.3cm]{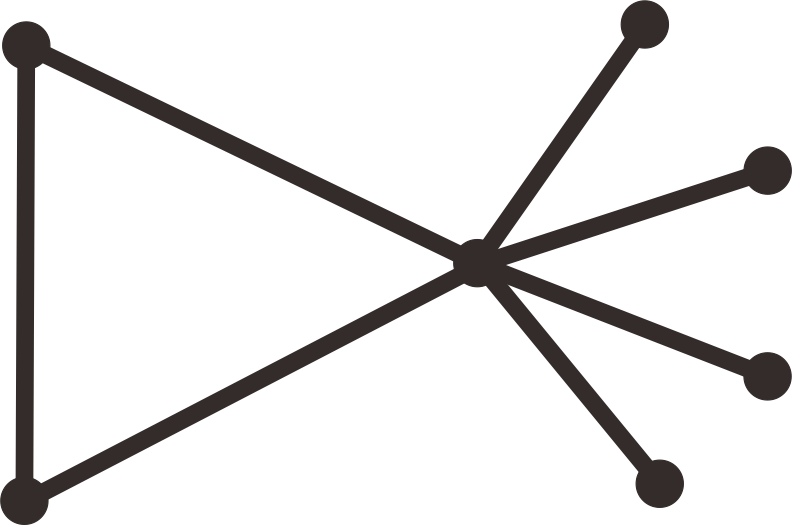}}
	\end{minipage}
    & 3.4526
    &
    &
    \\[2.5pt]
     \hline
  \end{tabular}
  \label{table 5}
\end{table}

\end{document}